\documentclass[10pt,a4paper]{amsart}
\RequirePackage[OT1]{fontenc}
\RequirePackage{amsthm,amsmath,amssymb,mathrsfs, enumerate, tikz, textcomp}
\RequirePackage[colorlinks,citecolor=blue,urlcolor=blue]{hyperref}
\RequirePackage{hyphenat, fullpage}
\newtheorem{theorem}{Theorem}[section]
\newtheorem{corollary}[theorem]{Corollary}
\newtheorem{lemma}[theorem]{Lemma}
\newtheorem{prop}[theorem]{Proposition}

\newtheorem{definition}[theorem]{Definition}
\newtheorem{example}[theorem]{Example}
\newtheorem{probl}[theorem]{Problem}
\newtheorem{remark}[theorem]{Remark}
\newcommand{\Ii}{1\!\!1}
\newcommand{\be}{\begin{equation}}
\newcommand{\ee}{\end{equation}}
\usepackage{graphicx}
\def\RR{\mathbb{R}}
\def\NN{\mathbb{N}}
\def\r{{\bf r}}\def\r{{\bf r}}
\def\y{{\bf y}}\def\z{{\bf z}}\def\x{{\bf x}}
\def\C{{\mathcal C}}
\def\F{{\mathcal F}}
\def\M{{\mathcal M}}
\def\R{{\mathcal R}}

\usepackage{xcolor}
\usepackage{ulem}

\newcommand\copyrighttext{%
  \footnotesize 
\textcopyright 2019. This manuscript version is made available under the CC-BY-NC-ND~4.0 license \href{http://creativecommons.org/licenses/by-nc-nd/4.0/}{http://creativecommons.org/licenses/by-nc-nd/4.0/} \\ \vspace{-0.2cm}

The final publication is available at \href{https://doi.org/10.1016/j.jmaa.2019.123551}{https://doi.org/10.1016/j.jmaa.2019.123551} 
}
\newcommand\copyrightnotice{%
\begin{tikzpicture}
\node[text width=15cm, anchor=south] at (current page.south) {\fbox{\parbox{\dimexpr\textwidth-\fboxsep-\fboxrule\relax}{\copyrighttext}}};
\end{tikzpicture}%
}

\title[The KMP on subsets of probabilities and point configurations]{The full moment problem on subsets of probabilities and point configurations.}
\author[M.Infusino and T. Kuna]{Maria Infusino and Tobias Kuna}
\address[Maria Infusino]{Fachbereich Mathematik und Statistik, Universit\"at
Konstanz, Universit\"atsstra{\ss}e 10, Konstanz 78457, Germany.}
\email{maria.infusino@uni-konstanz.de}
\address[Tobias Kuna]{Department of Mathematics and Statistics, University of Reading,
Whiteknights, PO Box 220, Reading RG6~6AX, United Kingdom}
\email{t.kuna@reading.ac.uk}

\keywords{Infinite dimensional moment problem, semi-algebraic set, random measure, point process, correlation.\\ \vspace{0.3cm} \copyrightnotice}

\begin{document}

\maketitle\vspace{-1cm}
\begin{abstract}
The aim of this paper is to study the full $K-$moment problem for measures supported on some particular non-linear subsets $K$ of an infinite dimensional vector space. We focus on the case of random measures, that is $K$ is a subset of all non-negative Radon measures on $\mathbb{R}^d$. We consider as $K$ the space of sub-probabilities, probabilities and point configurations on $\mathbb{R}^d$. For each of these spaces we provide at least one representation as a generalized basic closed semi-algebraic set to apply the main result in [J. Funct. Anal., 267 (2014) no.5: 1382--1418]. We demonstrate that this main result can be significantly improved by further considerations based on the particular chosen representation of $K$. 
In the case when $K$ is a space of point configurations, the correlation functions (also known as factorial moment functions) are easier to handle than the ordinary moment functions. Hence, we additionally express the main results in terms of correlation functions.
\end{abstract}
\normalem

  \enlargethispage{0.6cm}

\section*{Introduction}

The classical full $d-$dimensional $K-$moment problem ($d-$KMP) asks whether a multi-sequence of real numbers $m=(m_{\alpha})_{\alpha\in\NN_0^d}$ is actually the moment sequence associated to a non-negative Radon measure whose support is contained in a fixed closed subset $K$ of $\RR^d$ (here $d\in\mathbb{N}$). 
If such a measure exists, it is said to be a $K-$representing measure for $m$.  
Already at an early stage of its history, the $d-$KMP was generalized to the case when $K$ is a subset of an infinite dimensional vector space.
We refer to this problem as the \emph{infinite dimensional $K-$moment problem} ($\infty-$KMP). A general treatment of the full $\infty-$KMP has been given e.g. in \cite{AJK15}, \cite[Chapter~5, Section~2]{BeKo88}, \cite{BS71}, \cite{Bor-Yng75}, \cite{GIKM}, \cite{GKM16}, \cite{Heg75}, {\cite{IKKM19}}, \cite[Section~12.5]{Schmu90}, {\cite{Schm18}}. Special infinite dimensional supports, particularly useful in applications, have been considered e.g. in \cite{IKR14}, \cite{Le75a}, \cite{Le75b}, (see also \cite{CaKuLeSp06}, \cite{Kor05}, \cite{KLS07}, \cite{KLS11}, \cite{MolLach15} for the truncated case, i.e. when the starting $m$ consists of finitely many elements).

In this article we focus on the \emph{full} $\infty-$KMP for some specific non-linear subsets $K$ of the space of non-negative Radon measures on $\RR^d$, namely the set of all sub-probabilities, probabilities and point configurations on $\RR^d$. In probabilistic language these KMPs correspond to the full moment problem for random sub-probabilities, probabilities, random point processes or point fields, respectively. The choice of $\RR^d$ as underlying space is just for simplicity and the results can be easily extended for example to differentiable manifolds.

The main distinguishing feature of this paper, in contrast to the aforementioned ones, is that specific properties of the considered $K$ are exploited to get solvability conditions for the KMP which are much weaker than the ones directly obtainable from general results. 

In order to have a joint convenient framework we consider all these $K$ as subsets of the vector space $\mathscr{D}'(\RR^d)$ of all generalized functions on $\RR^d$, which has favourable properties as dual of the space $\C_c^\infty(\RR^d)$ of all infinitely differentiable functions with compact support in~$\RR^d$ (for example, the $n-$th moment $m^{(n)}$ of a measure on $\mathscr{D}'(\RR^d)$ can be considered as a generalized function in $\mathscr{D}'(\RR^{dn})$). The starting point of our investigation is the observation that all the supports $K$ considered here can be determined by (infinitely many) polynomial constraints on $\mathscr{D}'(\RR^d)$. In the following we will refer to a subset of $\mathscr{D}'(\RR^d)$ having this property as a \emph{generalized basic closed semi-algebraic set} (g.b.c.s.s.). 
Clearly, such a subset $K$ of $\mathscr{D}'(\RR^d)$ can have more than one representation as g.b.c.s.s.. Polynomials on $\mathscr{D}'(\RR^d)$ are real valued functions which are reduced to usual polynomials when restricted to any finite dimensional subspace of $\mathscr{D}'(\RR^d)$. We denote by~$\mathscr{P}$ the set of all these polynomials on $\mathscr{D}'(\RR^d)$ such that any homogeneous polynomial $\eta \mapsto P(\eta)$ of degree $n$ can be written as the pairing of $\eta^{\otimes n}$ with an element from $\C_c^\infty(\RR^{dn})$. Hence, we do not only consider linear combinations of $n-$fold products of linear functionals on $\mathscr{D}'(\RR^d)$, but also their limits in a suitable topology and so the completion of the tensor algebra w.r.t. to such a topology. This gives us a much higher flexibility in the description of g.b.c.s.s..

Characterization of measures via moments are generally built up out of five different types of conditions:
\begin{enumerate}[I.]
\item\label{condpos} positivity conditions on the moment sequence;
\item\label{condquasi} conditions on the asymptotic behaviour of the moment sequence as the order goes to infinity;
\item\label{condsupp} conditions on the support $K$ of the representing measure;
\item\label{condreg} regularity properties of the moments as generalized functions;
\item\label{condgro} growth properties of the moments as generalized functions. 
\end{enumerate} 
Conditions of type \ref{condreg} and \ref{condgro} are only relevant for the infinite dimensional moment problem. The general aim in moment theory is to obtain characterizations of the solutions to a given KMP which are as weak as possible w.r.t.\!\! some combination of the above different types of conditions, since it seems unfeasible to get a result which is optimal in all types simultaneously.

In this paper, we exploit \ref{condsupp} to weaken the conditions of the other types. Our results are based on \cite[Theorem 2.3]{IKR14}, which is a general criterion to solve the $\infty-$KMP for g.b.c.s.s.\! $K$ of $\mathscr{D}'(\RR^d)$. Here we derive improved results exploiting the special structure of the $K$'s under consideration. The solvability conditions in \cite[Theorem~2.3]{IKR14} do not include either any condition of type \ref{condsupp} beside $K$ being a g.b.c.s.s.\! of $\mathscr{D}'(\RR^d)$ nor any conditions of type \ref{condreg} and \ref{condgro} except the $n-$th moment $m^{(n)}$ being a generalized function. However, the essential conditions required are the positive semi-definiteness of the putative moment sequence $(m^{(n)})_{n\in\NN_0}$ and of some shifted versions of it (i.e. conditions of type~\ref{condpos}) and a growth condition of Carleman type on certain bounds of the $m^{(n)}$'s (so a condition of type \ref{condquasi}). Note that the positivity conditions in \cite[Theorem 2.3]{IKR14} depend on the chosen representation of the $K$'s as g.b.c.s.s.. 
In the following we derive properties of the moments directly from a convenient representation of $K$ as g.b.c.s.s. without using the existence of a $K-$representing measure. In this way we can in some cases weaken the conditions of type \ref{condpos} and~\ref{condquasi} given in \cite[Theorem 2.3]{IKR14}, showing nevertheless the existence of a representing measure. Each of the cases considered here demonstrate the use of a different technique to improve the conditions of \cite[Theorem 2.3]{IKR14}. On the one hand these examples, interesting in its own right, show the power of the method developed in \cite{IKR14}. 
On the other hand, they pose further challenges and point at potential for further development beyond the general theory.  

Let us describe our results in more details following the structure of the paper.

In Section \ref{sec1}, we introduce some basic notions and formulate the KMP for $K\subseteq\mathscr{D}'(\RR^d)$. In Subsection~\ref{subsec2.1} we recall our previous result in \cite{IKR14}, which combines techniques from the finite and the infinite dimensional moment theory (see reference therein), in particular \cite{BeKo88}, \cite{BS71}, \cite{Las2011} and \cite{Schm91}. The case where $K$ is the space of sub-probabilities is treated in Subsection~\ref{ExSubProb}. Inspired by the results in \cite{Schm91}, we are able to {directly derive} from the positive semi-definiteness assumptions a bound on the sequence of moments{,} which guarantees its determinacy. Any measure supported on the sub-probabilities gives rise indeed to a determined moment sequence.  However, the catch here is that we establish this bound without using the existence of a representing measure. Hence, we get a theorem solving the moment problem for random sub-probabilities which essentially only involves conditions of type {\ref{condpos}, namely of positive semi-definite type, but not of type \ref{condquasi}}. This is possible for a carefully chosen representation of {the space of sub-probabilities} as a g.b.c.s.s.. To treat the case where $K$ is the space of probabilities, we show in Subsection~\ref{subsecprob} that it is necessary and sufficient to add a single extra condition containing only the moments of order zero and one. In Subsection~\ref{subsecpoint} we consider the space of point configurations as the set $K$. The polynomials we use to represent this {space} as a g.b.c.s.s.\! are well-known and give rise to 
the so-called \emph{factorial moment measures} in the theory of point processes, also known as \emph{correlation functions} in statistical mechanics. To the best of our knowledge, it was not known before that the non-negativity of polynomials of this class with non-negative coefficients characterizes the space of point configurations. We also study the case of simple point configurations, {i.e.} point configurations which have at most one point at the same position. As an extra condition we introduce here that the second moment function on the diagonal coincides with the first one. Surprisingly, in this case the number of conditions of type \ref{condpos} can be reduced to only considering shifts of the putative moment sequence by polynomials {of} degree {at most $2$ instead of all degrees}. Note that these positive semi-definiteness conditions are not sufficient to represent the space of simple point configurations as a g.b.c.s.s..

The determining condition in all these cases can be weakened to a Stieltjes type growth condition on the $m^{(n)}$'s which requires only a bound of the form $B \left(Cn^2 \ln(n)^2 \right)^n$ for some constant $B,C \geq 0$. This had been already proved by A.~Lenard~\cite{Le72} for the case of point configurations. Positive semi-definite type solvability conditions for the moment problem for point processes were given before only under determining conditions of the type $B \left(Cn \ln(n) \right)^n$ and additional positivity or regularity assumptions, cf. \cite{BeKoKuLy99} and \cite{KoKu99}, whereas the classical Ruelle bound corresponds to a growth of type $C^n$ in our notation. There exist models with strong clustering where the Stieltjes type condition is achieved, cf. \cite{KoKuSe08}.

The solvability conditions for KMP for point processes and point random fields are frequently formulated in terms of correlation functions instead of moment functions. Therefore, in Section~\ref{seccor} we rewrite our results of the previous sections in terms of correlation functions. To this aim we need to extend the ``harmonic analysis'' on point configuration spaces in \cite{KoKu99} to $\mathscr{D}'(\RR^d)$. Though this was already used in the past, we are not aware of a systematic exposition as in Subsection~\ref{secHAgen}. In Subsection~\ref{subseccor} we show that the sequence of putative correlation functions fulfills the Stieltjes determining condition if and only if the sequence of associated usual moments fulfills the Stieltjes determining condition {together with an extra} mild assumption. Moreover, one can show that a sufficient solvability condition is the positive semi-definiteness of a shifted sequence of correlation functions with some of its arguments fixed (see Theorem~\ref{Thm-Suff} vs \eqref{newcond*}). Certain technical proofs {and remarks} are moved to the appendix.

Let us briefly state some related results on the $\infty-$KMP for the special supports considered here, which are based on positive semi-definiteness. (Results for the $\infty-$KMP holding for generic $K$ have been mentioned in the first paragraph of this introduction.) The case of random measures, that is $K$ is the cone of all non-negative Radon measures, has been treated in \cite{S74} where the cone structure is used to improve the conditions of type \ref{condquasi}. The case when $K$ is the set of all sub-probabilities can be treated using the general result in \cite{GIKM} (see Appendix \ref{Rem-SP}). The moment problem for point processes and point random fields has a rich and long history starting with A.~Lenard showing an analogue of Riesz-Haviland's result in \cite{Le75a}. In the same period K.~Krickeberg \cite{Krick73} characterizes point processes via restrictions of moment functions to the diagonals. Beside the results in  \cite{BeKoKuLy99} and \cite{KoKu99} mentioned before, solutions to the moment problem on point configuration spaces using positive semi-definiteness have been formulated also in terms of the generating function of the correlation functions, the so-called Bogoliubov functional, see e.g. \cite{KoKuOl02} for a result under $L^1-$analyticity. The case of random discrete measures treated in \cite{KoKuLyt15} cannot be treated via g.b.c.s.s.\! as the required support is not even closed. 
All these works resolve in different ways the balancing among the conditions \ref{condpos} to \ref{condgro} but the solvability conditions they provide are not comparable with the ones in this article, in the sense that one cannot show the equivalence of the two sets of conditions without using that each of them guarantees the existence of a representing measure.

\section{Preliminaries}\label{sec1}

In this section we state the full $\infty-$KMP for $K$ closed subset of the space of all generalized functions on $\RR^d$ according to the notation used in \cite{IKR14}. 

Let us start by recalling some preliminary notations and definitions. For $Y\subseteq \RR^d$, we denote by $\mathcal{B}(Y)$ the Borel $\sigma$-algebra on $Y$, by $\C_c^\infty(Y)$ the space of all real valued infinitely differentiable functions on $\RR^d$ with compact support contained in $Y$ and by $\C_c^{+,\infty}(Y)$ the cone consisting of all non-negative functions in $\C_c^\infty(Y)$. We denote
by $1\!\!1_Y$ the indicator function for $Y$ and by $\NN_0$ the space of all non-negative integers. For any
$\r=(r_1,\ldots,r_d)\in\RR^d$ and $\alpha=(\alpha_1,\ldots,\alpha_d)\in\NN_0^d$ one defines $\r^\alpha:=r_1^{\alpha_1}\cdots r_d^{\alpha_d}$. Moreover, for any $\beta\in\NN_0^d$ the symbol $D^\beta$ denotes the weak partial derivative $\frac{\partial^{\left|\beta\right|}}{\partial r_1^{\beta_1}\cdots\partial r_d^{\beta_d}}$ where $\left|\beta\right|:=\sum_{i=1}^d\beta_i$.

{The classical topology considered on $\C_c^\infty(\RR^d)$ is the inductive topology $\tau_{ind}$, given by the standard construction of this space as the inductive limit of spaces of smooth functions with supports lying in an increasing sequence of compact subsets of $\RR^d$ (see e.g.~\cite[Chapter 13, Example II]{Tre67}, \cite[Definition~5.9]{IKR14}). We denote by $\mathscr{D}_{ind}(\RR^d)$ the space $\C_c^\infty(\RR^d)$ equipped with~$\tau_{ind}$. In this paper, we consider instead $\C_c^\infty(\RR^d)$ endowed with the projective topology~$\tau_{proj}$ defined as follows (see \cite[Chapter I, Section 3.10]{B86} for more details).}

{
\begin{definition}\label{TeoBerez}\ \\
Let $I$ be the set of all $k=(k_1, k_2(\r))$ such that $k_1\in\NN_0$, $k_2\in\C^\infty(\RR^d)$ with $k_2(\r)\geq 1$ for all $\r\in\RR^d$.  For each $k=(k_1, k_2(\r))\in I$, consider the weighted Sobolev space $W_2^{k}$ defined as the completion of $\C_c^\infty(\RR^d)$ w.r.t.\! the following weighted norm
$$
\|\varphi\|_{W_2^{k}}:=\left(\sum_{\stackrel{\beta\in\NN_0^d}{|\beta|\leq k_1}}\int_{\RR^d}\left |(D^{\beta}\varphi)(\r)\right |^2k_2(\r)d\r\right)^{\frac 12}.
$$
Then we define $$\mathscr{D}(\RR^d):=\projlim\limits_{k\in I}W_2^{k},
$$ and we denote by $\tau_{proj}$ the projective limit topology induced by this construction.
\end{definition}

As a set $\mathscr{D}(\RR^d)$ coincides with $\C_c^\infty(\RR^d)$ (see \cite[Chapter I, Theorem~3.9]{B86}).} The space $\mathscr{D}(\RR^d)$ is a locally convex vector space which is also \emph{nuclear} (see e.g. \cite[Chapter I, Theorem~3.9]{B86} for a proof of this result and \cite[Chapter~14, Sections 2.2--2.3]{BSU96} for more details about general nuclear spaces).  We denote by $\mathscr{D}'(\RR^d)$ the topological dual space of $\mathscr{D}(\RR^d)$ and by $\langle f, \eta\rangle$ the duality pairing between $\eta\in\mathscr{D}'(\RR^d)$ and $f\in\mathscr{D}(\RR^d)$ (see e.g.~\cite{B86, BeKo88, BSU96} for more details). When we assume more regularity on $\eta$, we tacitly extend the dual pairing to larger classes of test functions $f$ (e.g. if $\eta$ is a Radon measure on $\mathbb{R}^{d}$, then we consider $\langle f, \eta\rangle$ for any measurable function $f$ on $\mathbb{R}^{d}$). We equip $\mathscr{D}'(\RR^d)$ with the weak topology $\tau^{proj}_w$, that is, the smallest topology such that the mappings
$\eta \mapsto \langle f, \eta \rangle$ are continuous for all $f\in\mathscr{D}(\RR^d)$. 

Let us introduce now the main objects involved in the $K-$moment problem for subsets $K$ of $\mathscr{D}'(\RR^d)$. A generalized process on~$\mathscr{D}'(\RR^d)$ is a non-negative Radon measure $\mu$ defined on the Borel $\sigma-$algebra on~$\mathscr{D}'(\RR^d)$. Moreover, we say that a generalized process $\mu$ is \emph{concentrated on} a measurable subset $K\subseteq\mathscr{D}'(\RR^d)$ if $\mu\left(\mathscr{D}'(\RR^d)\setminus K\right)=0$.

\begin{definition}[$n-$th local moment]\ \\
Given $n\in\NN$, a generalized process $\mu$ on $\mathscr{D}'(\RR^d)$ has \emph{$n-$th local moment} (or local moment of order~$n$) if for every $f\in\C_c^\infty(\RR^d)$ we have
$$\int_{\mathscr{D}'(\RR^d)}|\langle f, \eta\rangle|^n \mu(d\eta)<\infty.$$
{If in addition the functional $f\mapsto\int_{\mathscr{D}'(\RR^d)}|\langle f, \eta\rangle|^n \mu(d\eta)$ is continuous on $\mathscr{D}(\RR^d)$, then we say that $\mu$ has \emph{continuous $n-$th local moment}.}
\end{definition}

{If a generalized process $\mu$ on $\mathscr{D}'(\RR^d)$ has continuous $n-$th local moment, then it is easy to show that
there exist a $k\in I$ and $C>0$ such that 
$$\left(\int_{\mathscr{D}'(\RR^d)}|\langle f, \eta\rangle|^n \mu(d\eta) \right)^{\frac 1n}\leq C\|f\|_{W_2^k}, \forall f\in \mathscr{D}(\RR^d).$$
By using H\"older's inequality, this in turn implies that 
$$\left|\int_{\mathscr{D}'(\RR^d)}\langle f_1, \eta\rangle\cdots\langle f_n, \eta\rangle \mu(d\eta) \right|\leq C^n\prod_{i=1}^n\|f_i\|_{W_2^k}, \forall f_1, \ldots, f_n\in \mathscr{D}(\RR^d).$$
This together with the nuclearity of $\mathscr{D}(\RR^d)$ allow us to apply the generalized version of the Kernel Theorem in \cite[Vol II, Chapter~14, Theorem~6.2]{BSU96} in order to get that there exist $j\in I$ and a unique symmetric $m^{(n)}_{\mu}\in \left(W_2^{-j}\right)^{\otimes n}$ (and so $m^{(n)}_{\mu}\in \mathscr{D}'(\RR^{dn})$) such that 
\be\label{dens}\langle f_1\otimes \cdots\otimes f_n, m^{(n)}_{\mu}\rangle =\int_{\mathscr{D}'(\RR^d)}\langle f_1, \eta\rangle\cdots\langle f_n, \eta\rangle \mu(d\eta),\ \forall f_1, \ldots, f_n\in \mathscr{D}(\RR^d).\ee

We also have that $m^{(n)}_{\mu}\in \mathscr{D}_{ind}'(\RR^{dn})$ as $\tau_{ind}$ is finer than $\tau_{proj}$. 
Moreover, by \cite[Corollary II.2.5]{Fernique} the map
$\begin{array}{lll}
 f^{(n)}&\mapsto&\int_{\mathscr{D}_{ind}'(\RR^{dn})}\langle f^{(n)}, \eta^{\otimes n}\rangle \mu( d\eta)
\end{array}$
is also in $\mathscr{D}_{ind}'(\RR^{dn})$. Since this map coincides with $m^{(n)}_{\mu}$ on a dense subset by \eqref{dens}, we have that
\begin{equation}\label{nMomTensor}
\langle f^{(n)}, m_{\mu}^{(n)} \rangle =\int_{\mathscr{D}'(\RR^d)} \langle f^{(n)}, \eta^{\otimes n} \rangle \mu(d\eta),\ \forall f^{(n)}\in\mathscr{D}(\RR^{dn}).
 \end{equation} 
 This justifies the following definition.
}
 
 \begin{definition}[$n-$th generalized moment function]\label{MomFun}\ \\
{Given $n\in\NN$ and a generalized process $\mu$ on $\mathscr{D}'(\RR^d)$ with continuous $n-$th local moment, its  \emph{$n-$th generalized moment function in the sense of $\mathscr{D}'(\RR^d)$} is the symmetric generalized function $m^{(n)}_{\mu}\in\mathscr{D}'(\RR^{dn})$ such that \eqref{nMomTensor} holds. By convention, $m_{\mu}^{(0)}:=\mu(\mathscr{D}'(\RR^d))$.}
\end{definition}

As described above to a generalized process $\mu$ can be associated the corresponding generalized moment functions given by \eqref{nMomTensor}. The moment problem, which in an infinite dimensional context is often called the \emph{realizability problem}, addresses exactly the inverse question.
\begin{probl}[Moment problem on $K\subseteq\mathscr{D}'(\RR^d)$ (or KMP)]\label{RealProb}\  \\
Let $K$ be a closed subset of $\mathscr{D}'(\RR^d)$, $N\in\mathbb{N}_0\cup\{+\infty\}$ and $m=(m^{(n)})_{n=0}^N$ such that each $m^{(n)}\in\mathscr{D}'(\RR^{dn})$ is a symmetric functional. Find a generalized process $\mu$ with generalized moment functions (in the sense of $\mathscr{D}'(\RR^d)$) of any order and concentrated on~$K$ such that $$m^{(n)}=m^{(n)}_\mu\quad\text{for}\,\,\, n=0,\ldots, N,$$
i.e. $m^{(n)}$ is the $n-$th generalized moment function of $\mu$ for $n=0,\ldots, N$.
\end{probl}
If such a measure $\mu$ does exist we say that $m=(m^{(n)})_{n=0}^N$ is \emph{realized} by $\mu$ on~$K$ or equivalently that $\mu$ is a \emph{$K-$representing} measure for the sequence $m$. Note that the definition requires that one finds a measure concentrated on~$K$ and not only on $\mathscr{D}'(\RR^d)$. In the case $N=\infty$ one speaks of the ``full KMP'', otherwise of the ``truncated KMP''. In the following we are going to focus only on the full case and so we address to it just as the moment problem.

To simplify the notation from now on we denote by
$\mathcal{M}^*(K)$ the collection of all {non-negative Radon measures} concentrated on a measurable subset $K$~of~$\mathscr{D}'(\RR^d)$ {(i.e.} generalized processes) with {continuous local moment} of any order and by $\mathcal{F}\left(\mathscr{D}'\right)$ the collection of all infinite sequences $(m^{(n)})_{n\in\NN_0}$ such that each $m^{(n)}\in\mathscr{D}'(\RR^{dn})$ is a symmetric functional of its $n$ variables\footnote{{The choice of the notation $\mathcal{F}\left(\mathscr{D}'\right)$ is motivated by the similarity with the classical symmetric Fock space.}}.

Let us introduce the version of the classical Riesz functional for Problem~\ref{RealProb}. Denote by $\mathscr{P}$ the set of all polynomials on $\mathscr{D}'(\RR^d)$ of the form 
$$
P(\eta) := \sum_{j=0}^N\langle p^{(j)},\eta^{\otimes j}\rangle,
$$
where $p^{(0)}\in\RR$ and $p^{(j)}\in\C^\infty_c(\RR^{dj})$, $j=1,\ldots,N$ with $N\in\NN$. 

\begin{definition}\label{DefFunct}\ \\
Given $m\in\mathcal{F}\left(\mathscr{D}'\right)$, we define its associated Riesz functional $L_m$ as 
\begin{eqnarray*}
L_m: &\mathscr{P}&\to\RR \\
 \  &P(\eta)=\sum\limits_{n=0}^N\langle p^{(n)} , \eta^{\otimes n}\rangle & \mapsto L_m(P):=\sum_{n=0}^N\langle p^{(n)} , m^{(n)}\rangle .
\end{eqnarray*}
\end{definition}
When $m$ is realized by $\mu\in\mathcal{M}^*(K)$, a direct calculation shows that for any $P\in\mathscr{P}$ we get 
$$
L_m(P)=\int_{K}{P(\eta)\,\mu(d\eta)}.
$$

Hence, an obvious property of type \ref{condpos} which is necessary for an element in $\mathcal{F}\left(\mathscr{D}'\right)$ to be the moment sequence of some measure on $\mathscr{D}'(\RR^d)$ is the following.
\begin{definition}[Positive semi-definite sequence]\label{PosSemiDef}\ \\
A sequence $\xi\in\mathcal{F}\left(\mathscr{D}'\right)$ is said to be \emph{positive semi-definite} if 
$$L_{\xi}(h^2)\geq 0, \forall h\in\mathscr{P}.$$
\end{definition}
This is a straightforward generalization of the classical notion of positive semi-definiteness of the Hankel matrices considered in the finite dimensional moment problem, that is equivalent to require that the associated Riesz functional is non-negative on squares of polynomials. 

\section{Realizability of Radon measures in terms of moment functions}
\subsection{Previous results}\label{subsec2.1}\ \\
In \cite[Theorem 2.3]{IKR14}, we derived necessary and sufficient conditions for the solvability of Problem~\ref{RealProb} in the full case when $K$ is a \emph{generalized basic closed semi-algebraic} (g.b.c.s.s.), namely \be\label{S-Semialg}
K=\bigcap_{i\in I}\left\{\eta\in\mathscr{D}'(\RR^d)| \ P_i(\eta)\geq 0\right\},
\ee 
where $I$ is an index set and $P_i\in\mathscr{P}$.  Note that the index set $I$ is not necessarily countable. When $I$ is finite, this definition agrees with the classical one of basic closed semi-algebraic subset.
Denote by $\mathscr{P}_{K}$ the set of all the polynomials $P_i$'s defining $K$. Then w.l.o.g. we can assume that~$0 \in I$ and that $P_0$ is the constant polynomial $P_0(\eta)=1$ for all $\eta\in\mathscr{D}'(\RR^d)$. \\ 

In this paper we are going to consider some well-known subsets of the infinite dimensional space $\mathcal{R}(\RR^d)$ of all non-negative Radon measures on~$\RR^d$ for which we will provide a representation as  g.b.c.s.s.. Recall that $\mathcal{R}(\RR^d)$ is the space of all non-negative Borel measures that are finite on compact subsets of $\RR^d$ and it is itself a g.b.c.s.s. of $\mathscr{D}'(\RR^d)$ (see \cite[Example 4.8]{IKR14}). As mentioned in \cite{IKR14}, using a result due to S.N.~\v{S}ifrin about the infinite dimensional moment problem on dual cones in nuclear spaces (see \cite{S74}), it is possible to obtain a version of \cite[Theorem~2.3]{IKR14} for the case when $K$ is a g.b.c.s.s.\! of $\mathcal{R}(\RR^d)$ (the latter is in fact the dual cone of $\C_c^{+,\infty}(\RR^d)$). Before stating this result, let us introduce a growth condition on a generic $\xi\in\mathcal{F}\left(\mathscr{D}'\right)$ which will turn out to be sufficient for the uniqueness of the representing measure for $\xi$ on such subsets of $\mathcal{R}(\RR^d)$.  

\begin{definition}[Stieltjes determining sequence]\label{DefSeq}\ \\ 
A sequence $\xi\in\mathcal{F}\left(\mathscr{D}'\right)$ is said to be \emph{Stieltjes determining} if and only if there exists a total subset $E$ of $\C_c^\infty(\RR^d)$ 
and a sequence $(\xi_n)_{n\in\NN_0}$ of real numbers such that
$$
\xi_0=\sqrt{|\xi^{(0)}|}\,\text{ and }\,\xi_n\geq\sqrt{\sup_{f_1,\ldots,f_{2n}\in E}|\langle f_1\otimes\cdots\otimes f_{2n},\xi^{(2n)}\rangle|},\, \forall\,n\geq 1.
$$
and the class $C\{\sqrt{\xi_n}\}$ is quasi-analytic.
\end{definition}
Note that the classical Stieltjes condition $\sum_{n=1}^\infty \xi_n^{-\frac{1}{2n}}=\infty$ is sufficient for the class $C\{\sqrt{\xi_n}\}$ being quasi-analytic. For discussions about the choice of $E$ see \cite[Lemma 4.5.]{IKR14}.

We are ready now to state the result mentioned above about the KMP for $K$ g.b.c.s.s. in $\R(\RR^d)$.

\begin{theorem}\label{MainThm1}\ \\
Let $m\in\mathcal{F}\left(\mathscr{D}'\right)$ be a Stieltjes determining sequence and $K\subseteq\R(\RR^d)$ be a generalized basic closed semi-algebraic set of the form~\eqref{S-Semialg}.
Then $m$ is realized by a unique $\mu\in\mathcal{M}^*(K)$ \underline{if and only if} the following hold
\begin{equation}\label{MainThmCond1-1}
L_m(h^2)\geq 0, \,\,L_m(P_i h^2)\geq 0\,\,,\,\, \forall h\in\mathscr{P},\,\forall i\in I.
\end{equation}
\end{theorem}

Condition \eqref{MainThmCond1-1} is equivalent to require that the functional $L_m$ is non-negative on the quadratic module $\mathcal{Q}(\mathscr{P}_{K})$ associated to the representation \eqref{S-Semialg} of $K$, i.e.
$$\mathcal{Q}(\mathscr{P}_{K}):=\bigcup_{\stackrel{I_0\subset I}{|I_0|<\infty}}\left\{\sum_{i\in I_0} Q_i P_i\,:\, Q_i\in\Sigma\right\},$$
where $\Sigma$ is the set of all sum of squares of polynomials in $\mathscr{P}$.

{The connection between Theorem \ref{MainThm1} and \cite[Theorem 1.1]{GIKM} is discussed in Appendix~\ref{Rem-SP}.}

In this paper, we assume more regularity on the putative moment functions, that is, we require that they are all non-negative symmetric Radon measures, i.e. the starting sequence belongs to $\mathcal{F}\left(\mathcal{R}\right)$. This is actually the case in most of applications. One of the advantage of this additional assumption is that it allows us to rewrite the Stieltjes determinacy condition as follows.

\begin{definition}\label{WeighStieltj}\ \\
A sequence $\xi\in\mathcal{F}\left(\mathcal{R}\right)$ satisfies the \emph{weighted generalized Stieltjes condition} if for each $n\in\NN$ there exists a function $k^{(n)}_2\in\C^\infty(\RR^d)$ with $k^{(n)}_2(\r)\geq 1$ for all $\r\in\RR^d$ such that
\be\label{StieltjGen}
\sum_{n=1}^{\infty}\frac{1}{\sqrt{\sup\limits_{\stackrel{\z\in\RR^d}{\|\z\|\leq n}}\sup\limits_{\x \in [-1,1]^d}\sqrt{ \tilde{k}_2^{(n)}(\z+\x)}}\sqrt[4n]{\int_{\RR^{2nd}}\frac{\xi^{(2n)}(d\r_1,\ldots,d\r_{2n})}{\prod_{l=1}^{2n}k^{(2n)}_2(\r_l)}}}=\infty,
 \ee  
where $\tilde{k}^{(n)}_2\in\C^\infty(\RR^d)$ such that $\tilde{k}^{(n)}_2(\r)\geq\left|(D^\kappa k^{(n)}_2)(\r)\right|^2$ for all $|\kappa|\leq \lceil{\frac {d+1}{2}}\rceil$.
\end{definition}
As suggested by the name, the condition \eqref{StieltjGen} is an infinite-dimensional weighted version of the classical Stieltjes condition, which ensures the uniqueness of the solution to the one dimensional moment problem on $\RR^+$ (see \cite{Carl26, St1894}). Condition \eqref{StieltjGen} is sufficient for the sequence $\xi$ being Stietljes determining in the sense of Definition~\ref{DefSeq}. Using the weighted generalized Stieltjes condition, as already discussed in \cite[Section 4.2]{IKR14}, it is possible to prove the following corollary of Theorem~\ref{MainThm1}.
\begin{corollary}\label{MainThm2}\ \\
Let $m\in\mathcal{F}\left(\mathcal{R}\right)$ fulfill the weighted generalized Stieltjes condition in \eqref{StieltjGen} and let $K\subseteq\R(\RR^d)$ be a g.b.c.s.s.\! of the form~\eqref{S-Semialg}.
Then $m$ is realized by a unique $\mu\in\mathcal{M}^*(K)$ \underline{if and only if} the following hold
$$
L_m(h^2)\geq 0, \,\,L_m(P_i h^2)\geq 0,\,\,\, \forall h\in\mathscr{P},\,\forall i\in I,
$$
and for any $n\in\NN_0$ we have
\be\label{MainThmCond2}
\int_{\RR^{2nd}}\frac{m^{(2n)}(d\r_1,\ldots,d\r_{2n})}{\prod_{l=1}^{2n}k^{(2n)}_2(\r_l)}<\infty.
\ee
\end{corollary}

In the rest of this section, we are going to present different applications of Theorem~\ref{MainThm1} and Corollary~\ref{MainThm2}, giving concrete necessary and sufficient condition to solve the full Problem \ref{RealProb} for any starting sequence $m\in\mathcal{F}\left(\mathcal{R}\right)$ and for some well-known subset $K$ of $\R(\RR^d)$.

\subsection{The moment problem on the space of sub-probabilities $\mathcal{R}_{sub}(\RR^d)$}\label{ExSubProb}\ \\
We first provide a representation as g.b.c.s.s. of the set $\mathcal{R}_{sub}(\RR^d)$ of all sub-probabilities on $\RR^d$, which can be defined as 
$
\mathcal{R}_{sub}(\RR^d):=\{\eta\in\mathcal{R}(\RR^d): \eta(\RR^d)\leq 1\}.
$
From now on for any function $f\in\C_c(\RR^d)$ we denote by $\|f\|_\infty$ the supremum norm of $f$, i.e. $\|f\|_\infty:=\sup_{x\in\RR^d}|f(x)|$.
\begin{prop}\label{SubProbSemialg}\ \\
The set $\mathcal{R}_{sub}(\RR^d)$ is a g.b.c.s.s. of $\mathscr{D}'(\RR^d)$. More precisely, we get 
\footnotesize
\be\label{SubProbSemiAlg}
\mathcal{R}_{sub}(\RR^d)=\!\!\!\!\!\!\bigcap_{\psi\in\C_c^{+,\infty}(\RR^d)}\!\!\!\!\left\{\eta\in\mathscr{D}'(\RR^d): \Phi_\psi(\eta)\geq 0\right\}\,\,\, \cap\!\! \bigcap_{\stackrel{\varphi\in\C_c^{+,\infty}(\RR^d)}{\left\|\varphi\right\|_\infty\leq 1}}\!\!\!\!\left\{\eta\in\mathscr{D}'(\RR^d): \Upsilon_\varphi(\eta)\geq 0\right\},
\ee 
\normalsize
where $\Phi_\psi(\eta):=\langle \psi, \eta \rangle$ and $\Upsilon_\varphi(\eta):=1-\langle \varphi, \eta\rangle^2$. 
\end{prop}
\proof\ \\
Let us preliminarily recall that the set of all non-negative Radon measures on $\RR^d$ can be represented as a g.b.c.s.s. of $\mathscr{D}'(\RR^d)$ as follows (for a proof see \cite[Example 4.8]{IKR14}).
$$
\mathcal{R}(\RR^d)=\bigcap_{\psi\in\C_c^{+,\infty}(\RR^d)}\left\{\eta\in\mathscr{D}'(\RR^d): \Phi_\psi(\eta)\geq 0\right\}.
$$
Hence, the desired equality \eqref{SubProbSemiAlg} can be simply rewritten as
\be\label{SubProbSemiAlgI}
\mathcal{R}_{sub}(\RR^d)=\mathcal{R}(\RR^d) \,\,\cap\!\!\!\! \bigcap_{\stackrel{\varphi\in\C_c^{+,\infty}(\RR^d)}{\left\|\varphi\right\|_\infty\leq 1}}\!\!\left\{\eta\in\mathscr{D}'(\RR^d): \Upsilon_\varphi(\eta)\geq 0\right\}.
\ee 
Let $\eta\in\mathcal{R}_{sub}(\RR^d)$. Then for any $\varphi\in\C_c^{+,\infty}(\RR^d)$ with $\left\|\varphi\right\|_\infty\leq 1$ we have
$0\leq \langle \varphi,\eta\rangle\leq 1$,
which implies $\Upsilon_\varphi(\eta)\geq 0$.

Conversely, let $\eta$ be an element of the right-hand side of \eqref{SubProbSemiAlgI}. Then clearly $\eta\in\R(\RR^d)$ and for any $\varphi\in\C_c^{+,\infty}(\RR^d)$ with $\left\|\varphi\right\|_\infty\leq 1$ we easily get that
\be\label{eqSub}
0\leq \langle \varphi,\eta\rangle\leq 1.
\ee To prove $\eta\in\mathcal{R}_{sub}(\RR^d)$, it remains to show that $\eta(\RR^d)=\langle \Ii_{\RR^d},\eta\rangle\leq 1$. \\
For any positive real number $R$ let us define a function $\chi_R$ such that
\be\label{ChiR}
\chi_R\in\C_c^{+,\infty}(\RR^d)\, \text{ and }\, \chi_R(\r):=\left\{\begin{array}{ll}
1 & \text{if } |\r|\leq R\\
0 & \text{if } |\r|\geq R+1.
\end{array}
\right. 
\ee
Note that the function $\Ii_{\RR^d}$ can be approximated pointwise by the increasing sequence of functions $\{\chi_R\}_{R\in\RR^+}\subset\C_c^{+,\infty}(\RR^d)$ whose elements are s.t. $\left\|\chi_R\right\|_\infty=1$. Hence, by using the monotone convergence theorem and \eqref{eqSub}, we have that
$$\eta(\RR^d)=\langle \Ii_{\RR^d},\eta \rangle=\langle \lim_{R\to\infty} \chi_R,\eta \rangle= \lim_{R\to\infty}\langle \chi_R,\eta \rangle \leq 1.$$
\endproof

{Applying Corollary~\ref{MainThm2} for the case $K=\mathcal{R}_{sub}(\RR^d)$ and exploiting the representation \eqref{SubProbSemiAlg}, we are able to drop the conditions~\eqref{StieltjGen} and~\eqref{MainThmCond2} to obtain the following result.}

\begin{theorem}\label{ThmSP}\ \\
A sequence $m\in\mathcal{F}\left(\mathcal{R}\right)$ is realized by a unique $\mu\in\mathcal{M}^*(\mathcal{R}_{sub}(\RR^d))$  \underline{if and only if} the following inequalities hold
\begin{eqnarray}
&\ & L_m(h^2)\geq 0,\,\, \forall h\in\mathscr{P},\label{(a)}\\
&\ & L_m(\Phi_\psi h^2)\geq 0,\,\,\forall h\in\mathscr{P},\, \forall \psi\in\C_c^{+,\infty}(\RR^d),\label{(b)}\\
&\ & L_m(\Upsilon_\varphi h^2)\geq 0,\,\,\forall h\in\mathscr{P},\, \forall \varphi\in\C_c^{+,\infty}(\RR^d)\,\text{with}\,\left\|\varphi\right\|_\infty\leq 1,\label{(c)}
\end{eqnarray}
where $\Phi_\psi(\eta):=\langle \psi, \eta \rangle$ and $\Upsilon_\varphi(\eta):=1-\langle \varphi, \eta\rangle^2$. 
\end{theorem}
\proof\ \\
\noindent\textbf{Sufficiency}\ \\
Assume that \eqref{(a)}, \eqref{(b)} and \eqref{(c)} are fulfilled and let us show that \eqref{StieltjGen} and \eqref{MainThmCond2} hold for the function $k_2^{(n)}\equiv1,\ \forall n\in\NN.$
In fact, for any $n\in\NN$ and for any $\varphi\in\C_c^{+,\infty}(\RR^d)$ with $\left\|\varphi\right\|_\infty\leq 1$, we can apply \eqref{(c)} to $h(\eta)=\langle \varphi, \eta\rangle^{(n-1)}$. This implies that $$
L_m(\langle \varphi, \eta\rangle^{2n})=L_m(\langle \varphi, \eta\rangle^{2(n-1)}\langle \varphi, \eta\rangle^{2})\leq L_m(\langle \varphi, \eta\rangle^{2(n-1)}),
$$
and iterating, we get that
$$
L_m(\langle \varphi, \eta\rangle^{2n})\leq L_m(1).
$$
Consequently, for any real positive constant $R$, if we take in the previous inequality $\varphi=\chi_R$ as defined in \eqref{ChiR} then we have that 
$$\int_{\RR^{2nd}}\prod_{i=1}^{2n}\chi_R(\r_i){m^{(2n)}(d\r_1,\ldots,d\r_{2n})}=L_m(\langle \chi_R, \eta\rangle^{2n})\leq L_m(1).$$
Therefore, using the monotone convergence theorem as $R\to\infty$, we obtain
$$
\int_{\RR^{2nd}}{m^{(2n)}(d\r_1,\ldots,d\r_{2n})}\leq L_m(1)=m^{(0)}<\infty.
$$
Then, the conditions \eqref{MainThmCond2} and \eqref{StieltjGen} hold for $k_2^{(n)}\equiv1,\ \forall n\in\NN.$ Hence, by Corollary~\ref{MainThm2} the sequence $m$ is realized by a unique $\mu\in\mathcal{F}(\mathcal{R}_{sub}(\RR^d))$.

\noindent\textbf{Necessity}\ \\
The necessity of \eqref{(a)}, \eqref{(b)} and \eqref{(c)} follows by the simple observation that integrals of non-negative functions w.r.t.\! a non-negative measure are always non-negative. \\
\endproof

This proof was inspired by the results in \cite{Schm91} about the moment problem on a compact basic closed semi-algebraic subset of $\RR^d$. In fact, the set $\mathcal{R}_{sub}(\RR^d)$ is a compact subset of $\mathcal{R}(\RR^d)$ w.r.t.\! the vague topology (see Corollary A2.6.V in \cite{DaVJ03}). However, the technique in \cite{Schm91} does not apply straightforwardly to the moment problem on $\mathcal{R}_{sub}(\RR^d)$, because it only applies to classical basic closed semi-algebraic sets (i.e. defined by finitely many polynomials), which is not a natural situation in the infinite dimensional case we are considering in this paper.

 {For a discussion about the relation between Theorem \ref{ThmSP} and \cite[Theorem 1.1]{GIKM} see Appendix~\ref{Rem-SP}.}

The representation \eqref{SubProbSemiAlg} is not unique. It is indeed possible to give other representations of $\mathcal{R}_{sub}(\RR^d)$ as g.b.c.s.s. using slight modifications in the proof of Proposition~\ref{SubProbSemialg}. For instance, we can get
$$
\mathcal{R}_{sub}(\RR^d)=\bigcap_{\stackrel{\varphi\in\C_c^{+,\infty}(\RR^d)}{\left\|\varphi\right\|_\infty\leq 1}}\left\{\eta\in\mathscr{D}'(\RR^d): \langle \varphi, \eta\rangle-\langle \varphi, \eta\rangle^2\geq 0\right\},
$$
or also
\be\label{AltRepres}
\mathcal{R}_{sub}(\RR^d)=\!\!\!\!\!\!\bigcap_{\varphi\in\C_c^{+,\infty}(\RR^d)}\!\!\!\!\!\!\!\left\{\eta\in\mathscr{D}'(\RR^d): \langle \varphi, \eta\rangle\geq 0\right\}\,\,\cap\!\!\!\!\bigcap_{\stackrel{\varphi\in\C_c^{+,\infty}(\RR^d)}{\left\|\varphi\right\|_\infty\leq 1}}\!\!\!\!\!\!\!\left\{\eta\in\mathscr{D}'(\RR^d): 1-\langle \varphi, \eta\rangle\geq 0\right\}.
\ee
\normalsize
Depending on the choice of the representation, we get different versions of Corollary~\ref{MainThm2} for $\mathcal{R}_{sub}(\RR^d)$.

For instance, using the representation \eqref{AltRepres} in Corollary~\ref{MainThm2}, we obtain:
\begin{corollary}\ \\
Let $m\in\mathcal{F}\left(\mathcal{R}\right)$ fulfill \eqref{StieltjGen} and let $\mathcal{R}_{sub}(\RR^d)$ be represented as in \eqref{AltRepres}. Then $m$ is realized by a unique $\mu\in\mathcal{M}^*(\mathcal{R}_{sub}(\RR^d))$ if and only if \eqref{(a)}, \eqref{(b)}, \eqref{MainThmCond2} and the following hold:
\be\label{CondizioneReplace}
\hspace{-0.3cm}L_m( \Theta_\varphi h^2)\geq 0,\,\,\forall h\in\mathscr{P},\, \forall \varphi\in\C_c^{+,\infty}(\RR^d)\,\text{with}\,\left\|\varphi\right\|_\infty\leq 1,\ee
where $\Theta_\varphi(\eta):=1-\langle\varphi,\eta\rangle.$
\end{corollary}

Note that here we did not manage to drop \eqref{StieltjGen} and \eqref{MainThmCond2}, because the trick used in the proof of Theorem~\ref{ThmSP} does not work for the representation \eqref{AltRepres}.\\

The conditions \eqref{(a)}, \eqref{(b)} and \eqref{CondizioneReplace} can be rewritten more explicitly in terms of moment measures as 
$$\sum_{i,j}\langle h^{(i)}\otimes h^{(j)},\, m^{(i+j)}\rangle\geq 0,$$
$$\sum_{i,j}\langle h^{(i)}\otimes h^{(j)}\otimes\psi,\, m^{(i+j+1)}\rangle\geq 0,$$
$$\sum_{i,j}\langle h^{(i)}\otimes h^{(j)},\, m^{(i+j)}\rangle-\sum_{i,j}\langle h^{(i)}\otimes h^{(j)}\otimes\varphi,\, m^{(i+j+1)}\rangle\geq 0,$$
for all  $h^{(i)}\in\C_c^\infty(\RR^{id})$, $\psi\in\C_c^{+,\infty}(\RR^d)$ and $\varphi\in\C_c^{+,\infty}(\RR^d)$ with $\left\|\varphi\right\|_\infty\leq 1$. 

In particular, if each $m^{(n)}$ has a density $\alpha^{(n)}\in L^1(\RR^{dn},\lambda)$ w.r.t the Lebesgue measure $\lambda$ on $\RR^{dn}$, then \eqref{(a)}, \eqref{(b)} and \eqref{CondizioneReplace} respectively mean that $(\alpha^{(n)})_{n\in\NN_0}$ is positive semi-definite and for $\lambda-$almost all $\y\in\RR^d$ the sequence $(\alpha^{(n+1)}(\cdot,\y))_{n\in\NN_0}$ and $(\alpha^{(n)}(\cdot)-\alpha^{(n+1)}(\cdot,\y))_{n\in\NN_0}$ are positive semi-definite. 

This reformulation makes clear the analogy with the Hausdorff moment problem as treated in \cite{Devi53}, where $[0,1]$ is represented like
$$[0,1]=\{x\in\RR: x\geq 0\}\cap \{x\in\RR: 1-x\geq 0\}$$ and so necessary and sufficient conditions to solve the $[0,1]-$moment problem for a sequence of reals $(m_n)_{n\in\NN_0}$ are that $(m_n)_{n\in\NN_0}$, $(m_{n+1})_{n\in\NN_0}$ and $(m_n-m_{n+1})_{n\in\NN_0}$ are positive semi-definite. Also in this case we get different conditions on $(m_n)_{n\in\NN_0}$ depending on the representation we choose for $[0,1]$ as basic closed semi-algebraic subset of the real line (see \cite{BeMa84}).

\subsection{The moment problem on the space of probabilities $\mathcal{R}_{prob}(\RR^d)$}\label{subsecprob}\ \\
{The set $\mathcal{R}_{prob}(\RR^d)$ of all probabilities, i.e. $\mathcal{R}_{prob}(\RR^d):=\{\eta\in\mathcal{R}(\RR^d): \eta(\RR^d)= 1\}$, can be also represented as a g.b.c.s.s.\! of $\mathscr{D}'(\RR^d)$ defined by three infinite families of polynomials in $\mathscr{P}$. Hence, one can apply Corollary~\ref{MainThm2} for $K=\mathcal{R}_{prob}(\RR^d)$. In this subsection we instead treat $\mathcal{R}_{prob}(\RR^d)$ as a subset of $\mathcal{R}_{sub}(\RR^d)$ and apply Theorem~\ref{ThmSP}. This together with an additional trick brings the advantage that we can replace the infinitely many conditions coming from the third family mentioned above with a single condition~\eqref{In5}.}
\begin{theorem}\ \\
A sequence $m\in\mathcal{F}\left(\mathcal{R}\right)$ is realized by a unique $\mu\in\mathcal{M}^*(\mathcal{R}_{prob}(\RR^d))$ \underline{if and only if} the following inequalities hold
\begin{eqnarray}
&\ & L_m(h^2)\geq 0\,\,,\,\, \forall h\in\mathscr{P},\label{In1}\\
&\ & L_m(\Phi_\psi h^2)\geq 0\,\,,\forall h\in\mathscr{P},\,\, \forall \psi\in\C_c^{+,\infty}(\RR^d),\label{In2}\\
&\ & L_m(\Upsilon_\varphi h^2)\geq 0\,\,,\forall h\in\mathscr{P},\,\, \forall \varphi\in\C_c^{+,\infty}(\RR^d)\,\text{with}\,\left\|\varphi\right\|_\infty\leq 1,\label{In3}\\
&\ & m^{(1)}(\RR^d)=m^{(0)},\label{In5}
\end{eqnarray}
\normalsize
where $\Phi_\psi(\eta):=\langle \psi, \eta \rangle$, $\Upsilon_\varphi(\eta):=1-\langle \varphi, \eta\rangle^2$. 
\end{theorem}
\proof\ \\
\textbf{Necessity}\ \\
W.l.o.g. we can assume that the sequence $m$ is realized by a probability $\mu$ concentrated on~$\mathcal{R}_{prob}(\RR^d)$. This means that $m^{(0)}=1$. Moreover, $m$ is also realized on $\mathcal{R}_{sub}(\RR^d)\supset\mathcal{R}_{prob}(\RR^d)$ by $\mu$. Hence, Theorem~\ref{ThmSP} implies that \eqref{In1}, \eqref{In2} and \eqref{In3} hold. The condition \eqref{In5} easily follows from the assumption of the realizability of $m$ by approximating $\Ii_{\RR^d}$ with the increasing sequence of functions $\{\chi_R\}_{R\in\RR^+}\subset\C_c^{+,\infty}(\RR^d)$ defined in \eqref{ChiR}.\\
\noindent\textbf{Sufficiency}\ \\
Due to Theorem~\ref{ThmSP}, the assumptions \eqref{In1}, \eqref{In2} and \eqref{In3} imply that there exists a unique $\mu\in\mathcal{M}^*(\mathcal{R}_{sub}(\RR^d))$ representing $m$. Since $\mu$ is finite, we can assume w.l.o.g. that it is a probability on $\mathcal{R}_{sub}(\RR^d)$. It remains to prove that actually $\mu$ is concentrated on $\mathcal{R}_{prob}(\RR^d)$, i.e.
\be\label{suppOnProb}
\mu(\mathcal{R}_{prob}(\RR^d))=1.
\ee
On the one hand, as $\eta\in\mathcal{R}_{sub}(\RR^d)$, we get \be\label{PosIntegrand}
1-\langle \Ii_{\RR^d},\eta\rangle\geq 0.
\ee
On the other hand, approximating pointwisely the function $\Ii_{\RR^d}$ by an increasing sequence of functions $\{\chi_R\}_{R\in\RR^+}\subset\C_c^{+,\infty}(\RR^d)$ with $\left\|\chi_R\right\|_\infty=1$ (see~\eqref{ChiR} for the definition of $\chi_R$) and using the monotone convergence theorem together with the fact that $\mu$ is an $\mathcal{R}_{sub}(\RR^d)-$representing measure for $m$, we have that
$$\int_{\mathcal{R}_{sub}(\RR^d)}(1-\langle \Ii_{\RR^d},\eta\rangle)\mu(d\eta)=1-m^{(1)}(\RR^d)=0,$$
where in the last equality we used~\eqref{In5}.
Since $\mu$ is non-negative and by \eqref{PosIntegrand} the integrand is also non-negative on $\mathcal{R}_{sub}(\RR^d)$, the previous equation implies that $\mu-$a.s. in $\mathcal{R}_{sub}(\RR^d)$ we have
$1-\langle \Ii_{\RR^d},\eta\rangle=0$,
which is equivalent to \eqref{suppOnProb}.\\
\endproof

\subsection{The moment problem on point configuration spaces}\label{subsecpoint}\ \\
Let us preliminarily give a brief introduction to point configuration spaces (see~\cite{KoKu99}). For any subset $Y\in\mathcal{B}(\RR^d)$ and for any $n\in \NN_0$, we
define the\textit{\ space of multiple $n-$point configurations in $Y$} as 
$$
\ddot\Gamma^{(n)} _{0}(Y):=\{\delta_{x_1}+\cdots+\delta_{x_n} | x_i\in Y\},\,\forall n\in\NN$$
and $\ddot\Gamma _{0}^{(0)}(Y)$ as the set containing only the null-measure on $Y$.
To better understand the structure of $\ddot\Gamma^{(n)} _{0}(Y)$ we may use the following
natural mapping:
\begin{eqnarray}\label{symm}
\mathrm{sym}_Y^n:{Y^n} &\to& \ddot\Gamma _{0}^{(n)}(Y),\\
(x_1,\ldots ,x_n) &\mapsto&\delta_{x_1}+\cdots+\delta_{x_n}. \nonumber
\end{eqnarray}
Then it is clear that we can identify the space of multiple $n-$point 
configurations $\ddot\Gamma^{(n)} _{0}(Y)$ with the symmetrization of $Y^n$ w.r.t.\! the permutation group over $\{1,\ldots,n\}$ and endow it with the natural quotient topology.

We define the \textit{space of finite multiple configurations in $Y$} as $$
\ddot\Gamma _{0}(\RR^d):=\bigsqcup_{n\in \Bbb{N}_0}\ddot\Gamma _{0}^{(n)}(\RR^d)$$
equipped with the topology of disjoint union.

When we consider finite point configurations in $Y$ having in each site at most one point, we speak about \emph{finite simple point configurations in $Y$}. More precisely, for any $n\in \NN_0$, we
define the\textit{\ space of simple $n-$point configurations in $Y$} as $$
\Gamma^{(n)} _{0}(Y):=\{\eta\in \ddot\Gamma^{(n)} _{0}(Y)\, |\, \eta(\{x\})\leq 1,\ \forall\ x\in Y\},\,\forall n\in\NN
$$
and $\Gamma _{0}^{(0)}(Y)$ as the set containing only the null-measure on $Y$.
Therefore, using the mapping in \eqref{symm}, we can identify the space of simple $n-$point
configurations $\Gamma^{(n)} _{0}(Y)$ with the symmetrization of $\widetilde{Y^n}$ w.r.t.\! the permutation group over $\{1,\ldots,n\}$,
where 
\[
\widetilde{Y^n}:=\left\{ \left. (x_1,\ldots ,x_n)\in Y^n\right| \,x_k\neq
x_j\,\,\mathrm{if\,\,}k\neq j\right\} . 
\]
Using the following identification
$$
\sum_{i=1}^n\delta_{x_i}\longleftrightarrow \{x_1,\ldots, x_n\},
$$
we can also represent $\Gamma^{(n)} _{0}(Y)$ as a family of subsets of $Y$, that is,
\begin{equation}\label{fin-simpconf-set}\Gamma^{(n)} _{0}(Y)=\{\eta \subset Y\,|\,|\eta |=n\},\,\forall n\in \NN
\end{equation}
and $\Gamma _{0}^{(0)}(Y)=\{\emptyset\}$.

The \textit{space of finite simple configurations in $Y$} is then defined by
$$
\Gamma _{0}(Y):=\bigsqcup_{n\in \Bbb{N}_0}\Gamma _{0}^{(n)}(Y)$$
and can be equipped with the topology of disjoint union.\\

We are going to consider now the case of locally finite configurations of points in~$\RR^d$ and also in this case we will distinguish between multiple and simple configurations. Let us denote by $\mathcal{B}_c(\RR^d)$ the system of all sets in $\mathcal{B}(\RR^d)$ which are bounded and hence have compact closure. 

We define the \textit{space of multiple point configurations in $\RR^d$} as the set of all Radon measures on $\RR^d$ taking as values either a non-negative integer or infinity, i.e.
$$
\ddot\Gamma(\RR^d):=\left\{\eta\in\mathcal{R}(\RR^d)|\,\, \eta(B)\in\mathbb{N}_0, \forall B\in\mathcal{B}_c(\RR^d) \right\}.
$$
Any $\eta\in\ddot\Gamma(\RR^d)$ can be written as $\eta=\sum_{i\in I}\delta_{x_i}$ where $(x_i)_{i\in I}$  is such that $x_i\in\RR^d$ with $I$ either $\NN$ or a finite subset of $\NN$ and if $I=\NN$ then the sequence $(x_j)_{i\in I}$ has no accumulation points in $\RR^d$ (see \cite{DaVJ03}). 
This correspondence is one-to-one modulo relabelling of the points. The requirement that the sequence $(x_i)_{i\in I}$ has no accumulation points in $\RR^d$ corresponds to the condition that $\eta$ is a Radon measure on $\RR^d$. The space $\ddot\Gamma(\RR^d)$ is equipped with the \textit{vague topology} $\tau_v$, i.e. the weakest topology such that all the following functions are continuous
\begin{eqnarray*}
\ddot\Gamma(\RR^d)  &\rightarrow &\Bbb{R} \\
\eta  &\mapsto &\int f(x) \eta(dx),\quad\forall \ f\in \C_c(\RR^d).  \nonumber
\end{eqnarray*}

The \textit{space of simple point configurations in $\RR^d$} 
$$
\Gamma(\RR^d):=\{\eta\in\ddot\Gamma(\RR^d)|\,\forall\ \x\in\RR^d,\, \eta(\{ \x\})\in\{0,1\}\}.
$$
is considered with the relative vague topology induced by $(\ddot\Gamma(\RR^d), \tau_v)$.

From those definitions, it is then clear that point configurations in $\RR^d$ are subsets of $\mathscr{D}'(\RR^d)$ equipped with the weak topology.

There is also in this case a natural representation of $\Gamma (\RR^d)$ as a set of subsets in $\RR^d$, i.e.
$$
\Gamma(\RR^d) =\left\{ \left. \eta \subset \RR^d\right| \,|\eta \cap \Lambda|<\infty \,\,\forall\,\Lambda \in \mathcal{B}_c(\RR^d)\right\}.
$$
Indeed, any $\eta=\sum_{i\in I}\delta_{x_i}\in\ddot\Gamma(\RR^d)$ corresponds to $\{x_i\}_{i\in I}$.
One advantage of defining point configurations as Radon measures is the
ease of defining their powers, which is particularly convenient in the analysis of the moment problem on such spaces.

\subsubsection{The moment problem on the set of multiple point configurations $\ddot\Gamma(\RR^d)$}\ \\
The space of multiple point configurations of $\RR^d$ is also a g.b.c.s.s., namely we have the following representation.
\begin{prop}\label{AppN}\ \\
The set of multiple point configurations $\ddot\Gamma(\RR^d)$ on $\RR^d$
is a generalized basic closed semi-algebraic subset of $\mathscr{D}'(\RR^d)$. More precisely, we get that
\be\label{Nsemialg}
\ddot\Gamma(\RR^d)=\bigcap_{k\in\NN}\bigcap_{\varphi\in\C_c^{+,\infty}(\RR^{d})}\left\{\eta\in\mathscr{D}'(\RR^d): \, \langle \varphi^{\otimes k},\eta^{\odot k}\rangle\geq 0\right\}.
\ee
\end{prop}
The power $\eta^{\odot k}$ of a generalized function $\eta\in\mathscr{D}'(\RR^d)$ is called \emph{factorial power} and it is defined as follows. For any $f^{(n)}\in\C_c^{\infty}(\RR^{dn})$ and for any $n\in\NN$
\be\label{FactPower}
\langle f^{(n)} , \eta^{\odot n} \rangle := \sum_{k=1}^n \frac{(-1)^{n-k}}{k!} \sum_{\stackrel{n_1,\ldots, n_k \in \mathbb{N}}{n_1 + \ldots +n_k=n}}
\frac{n!}{n_1 \cdot \ldots \cdot n_k}
 \langle T_{n_1, \ldots, n_k}f^{(n)},
 \eta^{\otimes k} \rangle,
\ee
where
\be\label{FunctionT}
T_{n_1, \ldots, n_k}f^{(n)}(x_1, \ldots, x_k) := f^{(n)}(\underbrace{x_1,\ldots, x_1}_{n_1\text{ times}},\ldots, \underbrace{x_k,\ldots, x_k}_{n_k\text{ times}}).
\ee
The definition in \eqref{FactPower} shows that for any $n\in\NN$ and for any $f^{(n)}\in\C_c^{\infty}(\RR^{dn})$ the factorial power $\langle f^{(n)} , \eta^{\odot n} \rangle\in\mathscr{P}$. Note that when $f^{(n)}=f^{\otimes n}$ with $f\in \C^\infty_c(\RR^d)$ we have
\be\label{FactPowerFun}
T_{n_1, \ldots, n_k}f^{\otimes n}(x_1, \ldots, x_k) = f^{ n_1}(x_1)\cdots f^{ n_k}(x_k).
\ee
For example, in the cases $n=1$ and $n=2$ the previous definition gives 
$$ \langle f^{\otimes 1}, \eta^{\odot 1}\rangle=\langle f, \eta\rangle\quad\text{and}\quad 
\langle f^{\otimes 2}, \eta^{\odot 2}\rangle=\langle f, \eta\rangle^2- \langle f^2, \eta\rangle.$$
The name ``factorial power'' comes from the fact that for any $\eta\in\mathcal{R}(\RR^d)$ and for any measurable set $A$ 
 $$\langle\Ii_A^{\otimes n},\eta^{\odot n}\rangle=\eta(A)(\eta(A)-1)\cdots(\eta(A)-n+1).$$
 
Note that the definition of factorial power is very natural for point configurations in $\ddot\Gamma(\RR^d)$ (see \cite{KLS11}). In fact, using the representation $\eta=\sum_{i\in I}\delta_{x_i}$ for the elements in $\ddot\Gamma(\RR^d)$, \eqref{FactPower} becomes
\be\label{reprPointProcess}
\langle f^{\otimes n},\eta^{\odot n}\rangle=\sum'_{i_1,\ldots, i_n\in I} f(x_{i_1})\cdots f(x_{i_n}),
\ee
where $\sum'$ denotes a sum over distinct indices.

\proof(Proposition \ref{AppN})\ \\
Let $\eta\in\mathscr{D}'(\RR^d)$ such that for any $k\in\NN$ and for any $\varphi\in\C_c^{+,\infty}(\RR^{d})$
\be\label{Rel1}
\langle \varphi^{\otimes k},\eta^{\odot k}\rangle\geq 0.
\ee
In particular, the case $k=1$ implies that $\eta\in\mathcal{R}(\RR^d)$. Moreover, by a density argument, the condition \eqref{Rel1} also holds for $\varphi=\Ii_A$ with $A\in\mathcal{B}_c(\RR^d)$, i.e.
$$
0\leq\langle \Ii_A^{\otimes k},\eta^{\odot k}\rangle=\eta(A)(\eta(A)-1)\cdots(\eta(A)-k+1),\quad\forall k\in\NN,\forall A\in\mathcal{B}_c(\RR^d). 
$$
Hence, for any $A\in\mathcal{B}_c(\RR^d)$ we get that $\eta(A)\in\NN_0\cup\{+\infty\}$.\\
The other inclusion trivially follows from \eqref{reprPointProcess}.\\
\endproof
Using the representation \eqref{Nsemialg} and Theorem~\ref{MainThm1}, we have the following.
\begin{corollary}\ \\\label{ThmMultiplConfig}\
Let $m\in\mathcal{F}\left(\mathcal{R}\right)$ be a Stieltjes determining sequence. 
Then $m$ is realized by a unique $\mu\in\mathcal{M}^*(\ddot\Gamma(\RR^d))$ \underline{if and only if} the following hold:
\begin{eqnarray}
&\ &\hspace{-1.5cm} L_m(h^2)\geq 0,\,\, \forall h\in\mathscr{P},\nonumber\\
&\ &\hspace{-1.5cm} L_m(\Phi_{\varphi, k} h^2)\geq 0,\,\, \forall h\in\mathscr{P},\,\, \forall \varphi\in\C_c^{+,\infty}(\RR^d),\,\,\forall\, k\in\NN,\label{condCorr}
\end{eqnarray}
where $\Phi_{\varphi, k}(\eta):=\langle \varphi^{\otimes k},\eta^{\odot k}\rangle$.
\end{corollary}
\subsubsection{The moment problem on the set of simple point configurations $\Gamma(\RR^d)$}\ \\
The condition \eqref{condCorr} involves infinitely many polynomials of arbitrarily large degree. {However, we can show an analogue of Corollary \ref{ThmMultiplConfig} for the full $\Gamma(\RR^d)-$moment problem by requiring \eqref{condCorr} only for polynomials of at most second degree and by adding a single extra condition \eqref{condCorr4}.}

\begin{corollary}\label{ThmSimpleConfig}\ \\
Let $m\in\mathcal{F}\left(\mathcal{R}\right)$ be a Stieltjes determining sequence. 
Then $m$ is realized by a unique $\mu\in\mathcal{M}^*(\Gamma(\RR^d))$ \underline{if and only if} the following hold
\begin{eqnarray}
&\ & L_m(h^2)\geq 0,\,\, \forall h\in\mathscr{P}\label{condCorr0},\\
&\ & L_m(\Phi_{\varphi, 1} h^2)\geq 0,\,\,\forall h\in\mathscr{P},\,\, \forall \varphi\in\C_c^{+,\infty}(\RR^d),\label{condCorr1}\\
&\ & L_m(\Phi_{\varphi, 2} h^2)\geq 0,\,\,\forall h\in\mathscr{P},\,\, \forall \varphi\in\C_c^{+,\infty}(\RR^d),\label{condCorr2}\\
&\ & m^{(2)}(diag(\Lambda\times\Lambda))=m^{(1)}(\Lambda),\,\,\,\forall\,\Lambda\in\mathcal{B}(\RR^d)\,\text{compact},\label{condCorr4}
\end{eqnarray}
where $\Phi_{\varphi, 1}(\eta)=\langle \varphi,\eta\rangle$, $\Phi_{\varphi, 2}(\eta)=\langle \varphi^{\otimes 2},\eta^{\odot 2}\rangle=\langle \varphi,\eta\rangle^2-\langle \varphi^2,\eta\rangle$ and \\$diag(\Lambda\times\Lambda):=\{(x,x)\,|\, x\in\Lambda\}$.
\end{corollary}
\begin{remark}\label{prelRem}\ \\
By Theorem~\ref{MainThm1}, the conditions \eqref{condCorr0}, \eqref{condCorr1}, \eqref{condCorr2} are necessary and sufficient for the existence of a unique $\tilde{K}-$representing measure $\mu$ for the Stieltjes determining sequence $m$, where
$$\tilde{K}:=\!\!\!\!\bigcap_{\varphi\in\C_c^{+,\infty}(\RR^{d})}\!\!\!\!\left\{\eta\in\mathscr{D}'(\RR^d): \, \langle \varphi,\eta\rangle\geq 0\right\}\,\, \cap\!\!\!\! \bigcap_{\varphi\in\C_c^{+,\infty}(\RR^{d})}\!\!\!\!\left\{\eta\in\mathscr{D}'(\RR^d): \, \langle \varphi^{\otimes 2},\eta^{\odot 2}\rangle\geq 0\right\}.$$
Note that $\tilde{K}$ is not the required support in Corollary \ref{ThmSimpleConfig} as it is strictly larger than $\Gamma(\RR^d)$.
\end{remark}
\proof\ \\
\textbf{Sufficiency}\ \\
Let $\tilde{K}$ be as in Remark~\ref{prelRem}. By Theorem~\ref{MainThm1} there exists a unique $\tilde{K}-$representing measure $\mu$ for $m$. W.l.o.g.\! we can suppose that $\mu$ is a probability on $\tilde{K}$. Hence, it remains to show that $\mu$ is actually concentrated on $\Gamma(\RR^d)$. Let $\eta\in \tilde{K}$. Then for any $\varphi\in\C_c^{+,\infty}(\RR^d)$ we have:
$$\langle\varphi,\eta\rangle\geq 0\quad\text{and}\quad \langle \varphi^{\otimes 2},\eta^{\odot 2}\rangle\geq 0.$$
On the one hand, by a density argument, the previous conditions also hold for $\varphi=\Ii_A$ where $A\in\mathcal{B}(\RR^d)$ bounded, i.e.
$$
\left\{\begin{array}{l}
\eta(A)\geq 0 \\
\eta(A)(\eta(A)-1)\geq 0.
\end{array}\right.
$$
The latter relations imply that $\eta(A)\in\{0\}\cup[1,+\infty]$. Hence, for any $\eta\in \tilde{K}$ there exist $I\subseteq\NN$, $x_i\in\RR^d$ and real numbers $a_i\geq 1$ ($i\in I$) such that
\be\label{reprEta}
\eta=\sum_{i\in I}a_i\delta_{x_i},
\ee
where $I$ is either $\NN$ or a finite subset of $\NN$ and if $I=\NN$ then the sequence $(x_i)_{i\in I}$ has no accumulation points in $\RR^d$. \footnote{In fact, suppose that $supp(\eta)$ is not discrete then $\exists y\in supp(\eta)$ accumulation point. This means that there exists a compact neighbourhood $\Lambda$ of $y$ containing an infinite sequence $\{y_i\}_{i\in\NN}$ of points in $supp(\eta)$. Hence, $\eta(\Lambda)=\sum_{i=1}^\infty \eta(y_i)=\infty$ since $\eta(y_i)\geq 1$ ($\eta(y_i)$ cannot be zero because it $y_i$ in the support)}.

On the other hand, using \eqref{condCorr4}, the fact that $\mu$ is a $\tilde{K}$-representing measure for $m$ and that $\tilde{K}$ is a subset of Radon measures, we get via approximation arguments that for any $\Lambda\subset\RR^d$ measurable and compact 
$$0= m^{(2)}(diag(\Lambda\times\Lambda))-m^{(1)}(\Lambda)=\int_{\tilde{K}}(\langle \Ii_{diag(\Lambda\times\Lambda)},\eta^{\otimes 2}\rangle-\langle \Ii_\Lambda,\eta\rangle)\mu(d\eta).$$
As the integrand is non-negative on $\tilde{K}$, it follows that
$\langle \Ii_{diag(\Lambda\times\Lambda)},\eta^{\otimes 2}\rangle-\langle \Ii_\Lambda,\eta\rangle=0$ $\mu-$a.e. and so by~\eqref{reprEta} 
$$0=\sum_{\stackrel{i,j\in I}{x_i=x_j\in\Lambda}}a_ia_j-\sum_{\stackrel{i\in I}{x_i\in\Lambda}}a_i
=\sum_{\stackrel{i\in I}{x_i\in\Lambda}}a_i\left(\sum_{\stackrel{j\in I}{x_j=x_i}}a_j-1\right).$$
Since $a_i\geq1$ for all $i\in I$, we necessarily have that
$\sum\limits_{\stackrel{j\in I}{x_j=x_i}}a_j-1=0,$ namely
$$a_i-1+\sum_{\stackrel{j\neq i\in I}{x_j=x_i}}a_j=0.$$
The latter implies that $$\forall i\in I,\,\,a_i=1\quad\text{and}\quad \forall j,i\in I\,\,\text{with}\,\, j\neq i\,\text{we have }\,x_j\neq x_i.$$
Hence, we got that for $\mu-$ almost all $\eta\in K$
$$\eta=\sum_{i\in I}\delta_{x_i}\quad\text{and}\quad \eta(\{\x\})\in\{0,1\}$$
where $I$ is either $\NN$ or a finite subset of $\NN$ and if $I=\NN$ then the sequence $(x_i)_{i\in I}$ has no accumulation points in $\RR^d$.
This means that $\mu(\Gamma(\RR^d))=1.$

\noindent \textbf{Necessity}\ \\
By Remark~\ref{prelRem}, it only remains to show the condition \eqref{condCorr4}. Recall that for any $\eta\in\Gamma(\RR^d)$ there exist $I\subseteq\NN$ and $x_i\in\RR^d$ such that
$$\eta=\sum_{i\in I}\delta_{x_i}\quad\text{and}\quad \eta(\{\x\})\in\{0,1\}.$$
Therefore, for any $\Lambda\subset\RR^d$ measurable and compact 
$$
\langle \Ii_{diag(\Lambda\times\Lambda)},\eta^{\otimes 2}\rangle-\langle \Ii_\Lambda,\eta\rangle=\sum_{i,j\in I}\Ii_{diag(\Lambda\times\Lambda)}(x_i, x)-\sum_{i\in I} \Ii_\Lambda(x_i)=0.
$$
Hence, using approximation arguments and that $\mu$ is $\Gamma(\RR^d)-$representing for $m$, we get~\eqref{condCorr4}.\\
\endproof

\section{Realizability of Radon measures in terms of correlation functions}\label{seccor}
To simplify the notations in this section we will use the following abbreviations $\Gamma_0:=\Gamma_0(\RR^d)$, $\Gamma:=\Gamma(\RR^d)$, $\ddot\Gamma_0:=\ddot\Gamma_0(\RR^d)$, $\ddot\Gamma:=\ddot\Gamma(\RR^d)$.

\subsection{Harmonic analysis on generalized functions}\label{secHAgen}\ \\
As already mentioned in the introduction, we will need in the following some concepts from the so-called \emph{harmonic analysis on configuration spaces} developed in \cite{KoKu99}. However, since we aim to apply such notions to measures whose support is only a priori known to be contained in $\mathscr{D}'(\RR^d)$, we are going to provide a generalization of these concepts to our context.
\begin{definition}[$n-$th generalized correlation function]\label{corrFun}\ \\
Given $n\in\mathbb{N}$ and a generalized process $\mu$ on~$\mathscr
{D}'(\mathbb{R}^{d})$ with continuous $n-$th local moment, its \emph{$n-$th generalized correlation function in
the sense of~$\mathscr{D}'(\mathbb{R}^{d})$} is the symmetric
generalized function $\rho^{(n)}_{\mu}\in\mathscr{D}'(\mathbb{R}^{dn})$ such that

$$
\langle f^{(n)}, \rho_{\mu}^{(n)} \rangle=\int_{\mathscr
{D}'(\mathbb{R}^{d})} \langle f^{(n)}, \eta^{\odot n} \rangle\mu
(d\eta), \ \forall\ f^{(n)}\in{\mathcal C}_{c}^\infty(\mathbb{R}^{dn}).
$$
By convention, $\rho_{\mu}^{(0)}:=\mu(\mathscr{D}'(\mathbb{R}^{d}))$.
\end{definition}
In the previous definition the only change w.r.t.\! Definition \ref{MomFun} is that we consider a different basis for~$\mathscr{P}$, that is, we take $\frac{1}{n!}\langle f^{(n)}, \eta^{\odot n} \rangle$ instead of $\langle f^{(n)}, \eta^{\otimes n} \rangle$. This basis is different from the system of the Charlier polynomials commonly used in Poissonian analysis (see e.g.~\cite{CP90, IK88,KSSU98, NV95, P95} and \cite{KoKuOl02} for a detailed overview). 
Any polynomial in $\mathscr{P}$ can be written as:
$$P(\eta)=\frac{1}{j!}\sum_{j=0}^N\langle g^{(j)}, \eta^{\odot j} \rangle,$$
where $g^{(0)}\in\RR$ and $g^{(j)}\in\C_c^\infty(\RR^{dj})$, $j=1,\ldots, N$ with $N\in\NN$. W.l.o.g. each $g^{(j)}$ can be assumed to be a symmetric function of its $j$ variables in $\RR^d$. These symmetric coefficients are uniquely determined by $P$. One may introduce the following mapping which associates to a sequence of coefficients $(g^{(j)})_{j=0}^N$ the corresponding polynomial $P$. For convenience, we denote by $\F_f(\C_c^\infty)$ the collection of all sequences $G=(g^{(j)})_{j=0}^\infty$ where $g^{(j)}\in\C_c^\infty(\RR^{dj})$ is a symmetric function of its $j$ variables in $\RR^d$ and $g^{(j)}\neq 0$ only for finitely many $j$'s.

\begin{definition}[$K-$transform]\label{Ktransf}\ \\
For any $G=(g^{(j)})_{j=0}^\infty\in\F_f(\C_c^\infty)$, we define the $K-$transform of $G$ as the function $KG$ on $\mathscr{D}'(\RR^d)$ given by
$$
(KG)(\eta ):=\sum_{j=0}^\infty\frac{1}{j!}\langle g^{(j)}, \eta^{\odot j} \rangle,\, \forall\,\eta\in\mathscr{D}'(\RR^d).
$$
\end{definition}
The $K-$transform is well defined, because only finite many summands are different from 0. {One can give a more direct description of the $K-$transform whenever $\eta=\sum_{i\in I}\delta_{x_i}\in\ddot\Gamma$  (see Section \ref{subsecpoint}). Indeed, by using \eqref{reprPointProcess},  we have that
\begin{eqnarray*}
(KG)(\eta )\!\!\!\!&=&\!\!\!\!\sum_{j=0}^\infty\frac{1}{j!} \sum_{\stackrel{i_1,\ldots,i_j\in I}{i_1\neq\cdots\neq i_j}}g^{(j)}(x_{i_1}, \ldots , x_{i_j})
=\sum_{j=0}^\infty\sum_{\stackrel{i_1,\ldots,i_j\in I}{i_1<\cdots<i_j}}g^{(j)}(x_{i_1},\ldots,x_{i_j}).\\ \nonumber
\end{eqnarray*}
If $\eta\in\Gamma$ then one can also use the representation of $\Gamma$ as subsets given in \eqref{fin-simpconf-set} to obtain the following representation
\begin{eqnarray}\label{Ktransf-config}
(KG)(\eta )\!\!\!\!&=&\sum_{j=0}^\infty\sum_{\stackrel{\xi \subset \eta}{|\xi|=j}}\tilde{G}(\xi)=\sum_{\stackrel{\xi\subset\eta}{|\xi | < \infty }}\tilde{G}(\xi),
\end{eqnarray}
where $\tilde{G}:\Gamma_0\to\RR$ is defined as follows. Any $\xi\in\ddot\Gamma_0$ is of the form $\xi:=\{y_1,\ldots, y_n\}$ for some $n\in\NN$ and $y_i\in\RR^d$ (this representation is unique up to the relabelling of the $y_i$'s) then for such $\xi$ we define 
$$
\begin{array}{cccl}
\tilde{G}:&\Gamma_0&\to&\RR \\
 \ &\xi=\{y_1,\ldots, y_n\} & \mapsto & \tilde{G}(\xi):=g^{(n)}(y_1,\ldots,y_n)\\
\end{array}
$$
By convention, $\tilde{G}(\emptyset)=g^{(0)}$. The expression in \eqref{Ktransf-config} is what is used in \cite{KoKu99} as definition of $K-$transform.}

\begin{definition}[Convolution $\star$]\ \\
The convolution $\star:\F_f(\C_c^\infty)\times\F_f(\C_c^\infty)\to\F_f(\C_c^\infty)$ is defined as follows. For any $G, H\in\F_f(\C_c^\infty)$, we define $G \star H :=\left((G \star H)^{(j)}\right)_{j=0}^\infty\in \F_f(\C_c^\infty)$ to be such that for any $j\in\NN$
$$(G \star H)^{(j)}(x_J)=\sum_{\stackrel{J_1,J_2\subset J}{J_1\cup J_2=J}}G^{(|J_1|)}(x_{J_1})H^{(|J_2|)}(x_{J_2}),$$
where $J:=\{1,\ldots, j\}$ and, for any $I \subseteq J$, $x_I:=(x_i)_{i \in I}\in(\RR^{d})^{|I|}$. 
Note that $x_I$ is the equivalence class of all $|I|-$tuples with the same elements up to a relabelling.
\end{definition}

Let us prove here the analogous of Proposition 3.3 in \cite{KoKu99}.
\begin{prop}\label{PropConv}\ \\
Let $G,H\in \F_f(\C_c^\infty)$. Then 
$
K\left(G\star H\right) =KG\cdot KH.
$
\end{prop}
\proof\ \\
Due to the polarization identity, it suffices to consider the case $G=(g^{\otimes j})_{j=0}^\infty$ and $H=(h^{\otimes j})_{j=0}^\infty$ with $g,h\in\C_c^{\infty}(\RR^d)$. By using \eqref{FactPower} and \eqref{FactPowerFun} one easily get that for any $G=(g^{\otimes j})_{j=0}^\infty$ the following holds
$$KG(\eta)=\sum_{j=0}^\infty\frac{1}{j!}\langle g^{\otimes j}, \eta^{\odot j}\rangle=e^{\langle \ln(1+g), \eta\rangle}.$$

Then we have
\begin{eqnarray}\label{prodK}
KG(\eta)\cdot KH(\eta)&=&
e^{\langle \ln(1+g+h+gh), \eta\rangle}\\\nonumber
							  &=&\sum_{j=0}^\infty\frac{1}{j!}\langle (g+h+gh)^{\otimes j}, \eta^{\odot j}\rangle\\\nonumber
&=&\sum_{j=0}^\infty\frac{1}{j!}\sum_{j_1+j_2 + j_3 = j} \frac{j!}{j_1!j_2!j_3!} \langle g^{\otimes j_1}\hat{\otimes} (gh)^{\otimes j_2}\hat{\otimes} h^{\otimes j_3}, \eta^{\odot j}\rangle,
\end{eqnarray}
where $\hat{\otimes}$ denotes the symmetric tensor product.
Let us observe that 
\begin{eqnarray}\label{rewritingFunct}
& \ & \sum_{j_1+j_2 + j_3 = j} \frac{j!}{j_1!j_2!j_3!} g^{\otimes j_1}\hat{\otimes} (gh)^{\otimes j_2}\hat{\otimes} h^{\otimes j_3}(x_1,\ldots, x_j)\\\nonumber
&=&\sum_{j_1+j_2 + j_3 = j}\frac{1}{j_1!j_2!j_3!} \sum_{\pi\in S_j} \prod_{i=1}^{j_1+j_2}g(x_{\pi(i)}) \prod_{i=j_1+1}^{j}h(x_{\pi(i)})
\\ \nonumber
&=&\sum_{\stackrel{I_1,I_2\subset J}{I_1\cup I_2=J}}G^{(|I_1|)}(x_{I_1})H^{(|I_2|)}(x_{I_2})\\ \nonumber
&=& (G \star H)^{(j)}(x_1,\ldots, x_j)
\end{eqnarray}
where $S_j$ denotes the symmetric group on $\{1,\ldots, j\}$ and for any $\pi\in S_j$ we set  $J:=\{\pi(1),\ldots, \pi(j)\}$.

Then using \eqref{prodK} and \eqref{rewritingFunct} together with Definition~\ref{Ktransf}  we get the conclusion.
\endproof

Proposition \ref{PropConv} shows that $K$ is a $\RR-$algebra isomorphism between $(\F_f(\C_c^{\infty}), \star)$ and $(\mathscr{P}_{\C_c^\infty}(\mathscr{D}'), \cdot)$.\\

Given a sequence $\rho:=(\rho^{(n)})_{n\in\NN}\in\F(\mathscr{D'}(\RR^{d})$, we define the linear functional $\tilde{L}_{\rho}$ on $\mathscr{P}_{\C_c^\infty}(\mathscr{D}')$ by $\tilde{L}_{\rho}(\langle f^{(n)}, \eta^{\odot n}\rangle):=\langle f^{(n)}, \rho^{(n)}\rangle,$ $\forall\ n\in\NN_0,$ $\forall f^{(n)}\in\C_c^\infty(\RR^{dn})$, and $ \forall\eta\in\mathscr{D}'(\RR^d).$ If $\rho$ is realized by $\mu\in\M^*\left(\mathscr{D'}(\RR^{d})\right)$, then, by Definition \ref{corrFun}, we have
$$\tilde{L}_{\rho}(\langle f^{(n)}, \eta^{\odot n}\rangle)=\int_{\mathscr{D'}(\RR^{d})}\langle f^{(n)}, \eta^{\odot n}\rangle\mu(d\eta).$$

\begin{remark}\label{RemRightm}\ \\
To any $\rho:=(\rho^{(n)})_{n\in\NN}\in\F(\mathscr{D'}(\RR^{d})$ we can always associate the sequence $m=(m^{(n)})_{n\in\NN}\in\F(\mathscr{D'}(\RR^{d})$ such that $$\langle f^{(n)}, m^{(n)}\rangle:= \tilde{L}_{\rho}(\langle f^{(n)}, \eta^{\otimes n}\rangle).$$
In other words, $L_m=\tilde{L}_{\rho}$ on $\mathscr{P}_{\C_c^\infty}(\mathscr{D}')$.

For convenience, let us explicitly express $m$ in terms of $\rho$:
$$
\langle f^{(n)} , m^{(n)} \rangle = \sum_{k=1}^n \frac{1}{k!} \sum_{\stackrel{n_1,\ldots, n_k \in \mathbb{N}}{n_1 + \ldots +n_k=n}}
\frac{n!}{n_1! \cdot \ldots \cdot n_k!}
 \langle T_{n_1, \ldots, n_k}f^{(n)},
 \rho^{(k)} \rangle,
$$
where
$T_{n_1, \ldots, n_k}f^{(n)}$ is defined as in \eqref{FunctionT}.
\end{remark}

\subsection{The moment problem on the set of point configurations in $\RR^d$}\label{subseccor}\ \\
The positivity conditions in Corollaries \ref{ThmMultiplConfig} and \ref{ThmSimpleConfig} are of the form
\be\label{pdsCond}
L_m\left(\langle \varphi^{\otimes n}, \eta^{\odot n}\rangle h^2(\eta)\right)\geq 0,\,\, \forall h\in\mathscr{P}, \forall \varphi\in\C_c^{+,\infty}(\RR^d).\ee
We aim to express \eqref{pdsCond} in terms of correlation functions instead of moment functions. 
Given a sequence $\rho=(\rho^{(j)})_{j\in\NN}\in\F(\mathscr{D}'(\RR^d))$, we can always associate to it a sequence $m\in\F(\mathscr{D}'(\RR^d))$ as in Remark \ref{RemRightm}. Then for any $h\in\mathscr{P}$ and any $\varphi\in\C_c^{+,\infty}(\RR^d)$ we have
$$
L_m(\langle \varphi^{\otimes n}, \eta^{\odot n} \rangle h^2(\eta))=\tilde{L}_{\rho}(\langle \varphi^{\otimes n}, \eta^{\odot n} \rangle h^2(\eta)).
$$
It remains to consider under which conditions the Stieltjes determining property of a sequence $\rho$ (cf. Definition \ref{DefSeq}) is inherited by the associated $m$.
\begin{lemma}\label{detCorr}
Assume that a sequence $\rho=(\rho^{(j)})_{j\in\NN}\in\F(\mathscr{D}'(\RR^d))$ is Stieltjes determining. Then there exists a total subset $E$ of $\C_c^\infty(\RR^d)$ such that the class $C\{\sqrt[4]{\rho_{2n}}\}$ is quasi-analytic, where the sequence $(\rho_{n})_{n\in\NN_0}$ of real numbers is defined as
\begin{equation}\label{defCond}
\rho_0:=|\rho^{(0)}|\,\text{ and }\,\rho_n:= \sup_{f_1,\ldots,f_{n}\in E}|\langle f_1\otimes\cdots\otimes f_{n},\rho^{(n)}\rangle|,\, \forall\,n\geq 1.
\end{equation}
If the set $E$ is closed under multiplication and the sequence $\left(\frac{\rho_n}{n!}\right)_{n\in\NN}$ is almost-increasing (i.e. the exists $C>1$ s.t. $\frac{\rho_n}{n!}\leq C^{s}\left(\frac{\rho_s}{s!}\right)$ for any $n\leq s$), then the sequence $m$ associated to $\rho$ as in Remark \ref{RemRightm} is also Stieltjes determining.
\end{lemma}
The notion of almost-increasing sequence here introduced is weaker than the usual definition given e.g. in \cite{Mall59}.

\proof\ \\
The existence of $E$ and the quasi-analytic of the class $C\{\sqrt[4]{\rho_{2n}}\}$ directly follows from Definition \ref{DefSeq} applied to $\rho=(\rho^{(j)})_{j\in\NN}$.

Let $m=(m^{(j)})_{j\in\NN}$ be associated to $\rho$ as in Remark \ref{RemRightm} and let $$m_n:= \sqrt{\sup_{f_1,\ldots,f_{2n}\in E}|\langle f_1\otimes\cdots\otimes f_{2n},m^{(2n)}\rangle|},\, \forall\,n\geq 1.$$ W.l.o.g. we can assume that $m_n\geq 1$ for all $n\geq k$ for some $k\in\NN$ (otherwise it is clear that $\sum_{n=1}^\infty m_n^{-\frac{1}{2n}}=\infty$ holds and so that $m=(m^{(j)})_{j\in\NN}$ is Stieltjes determining). Then for any $g_1,\ldots, g_{2n}\in E$ and any $n\geq 1$ we have:
\begin{eqnarray*}
\sqrt{|\langle g_1\otimes\cdots\otimes g_{2n},m^{(2n)}\rangle|}&\!\!\!\leq\!\!\!&\sqrt{\sum_{k=1}^{2n}\frac{1}{k!}\!\!\sum_{\stackrel{n_1,\ldots, n_k \in \mathbb{N}}{n_1 + \ldots +n_k=2n}}
\!\!\!\frac{(2n)!}{n_1!\ldots n_k!}|\langle T_{n_1, \ldots, n_k}(g_1\otimes\cdots\otimes g_{2n}), \rho^{(k)}\rangle|}\\
&\leq&
\sqrt{\sum_{k=1}^{2n}\frac{1}{k!}\sum_{\stackrel{n_1,\ldots, n_k \in \mathbb{N}}{n_1 + \ldots +n_k=2n}}
\frac{(2n)!}{n_1!\ldots n_k!}\, \rho_k}
\leq\sqrt{\sum_{k=1}^{2n}k^{2n}\frac{\rho_k}{k!}}\\
&\leq&
\sqrt{(2n) (2n)^{2n}C^{2n}\frac{\rho_{2n}}{(2n)!}}
\leq\sqrt{\frac{(2eC)^{2n}}{\sqrt{\pi}}\rho_{2n}},
\end{eqnarray*}
where the constant $C$ is the one appearing in the almost increasing assumption and in the last step we have made use of the Stirling formula.

Hence, we get $m_n\leq\sqrt{\frac{(2eC)^{2n}}{\sqrt{\pi}} \rho_{2n}}$ and so the quasi-analyticity of the class $C\{\sqrt[4]{\rho_{2n}}\}$ implies that $C\{\sqrt{m_{n}}\}$ is quasi-analytic, i.e. $m=(m^{(j)})_{j\in\NN}$ is Stieltjes determining.
\endproof

\begin{remark}
As we have already showed in \cite[Section 2]{IKR14}, the determining condition becomes very concrete whenever one can explicitly construct the set~$E$. The statements in \cite[Lemma~2.5 and Remark~2.6]{IKR14} hold when replacing each $m^{(n)}$ with~$\rho^{(n)}$. There we showed that a preferable choice for $E$ is when $\sup_{f\in E}\|f\|_{H_{k^{(n)}}}$ grows as little as possible and we provided an example for such an $E$. The difference in the context of correlations is that Lemma~\ref{detCorr} requires an $E$ which has the additional property to be closed under multiplication. Along the same lines of the proof of  \cite[Lemma~2.5]{IKR14}, we introduce below in Lemma \ref{lemmino} a concrete example of $E$ fulfilling all the requirements in Lemma \ref{detCorr} and such that $\sup_{f\in E}\|f\|_{H_{k^{(n)}}}$ has the same growth as in \cite[Lemma 2.5]{IKR14}. For convenience, we consider here only the case when $E\subset\mathscr{C}^\infty_c(\RR)$. The higher dimensional case follows straightforwardly.
\end{remark}

Let $(d_n)_{n\in\NN_0}$ be a positive sequence which is not quasi-analytic, then there exists a non-negative infinite differentiable function $\varphi$ with support $[-1,1]$ such that for all $x\in\RR$ and $n\in\NN_0$ holds $|\frac{d^n}{dx^n}  \varphi(x)| \leq d_n$ (see \cite{Ru74}). 

\begin{lemma}\label{lemmino}\ \\
Let $(d_n)_{n\in\NN_0}$ be a log-convex positive sequence which is not quasi-analytic, let $\varphi$ be as above. Define
\[
E_0 := \{ f_{y,p,j}(\cdot):=\varphi^j(\cdot -y) e^{ip\cdot}\ | \ y ,p\in \mathbb{Q}, \ j\in\NN \}.
\]
Then for any $y,p\in\mathbb{Q}$ and for any $j, n\in\NN_0$ we get
\[
 \| f_{y,p,j}\|_{H_{k^{(n)}}} \leq C_p^{k_1^{(n)}} d_{k_1^{(n)}} \sup_{x \in [-1,1]} \sqrt{ k_2^{(n)}(y+x)},
\]
where $C_p:= \sqrt{2}(1+|p|)$
and $E_0$ is total in $\mathscr{D}(\mathbb{R})$.
\end{lemma} 
\proof (see Appendix \ref{App1})\endproof

Combining Lemma \ref{detCorr} with Corollary~\ref{ThmMultiplConfig} we get:
\begin{theorem}\label{ThmMultiplConfig-Corr}\ \\
Let $\rho\in\mathcal{F}\left(\mathcal{R}\right)$ be a Stieltjes determining sequence s.t. the corresponding $\left(\frac{\rho_n}{n!}\right)_{n\in\NN}$ is almost-increasing (see \eqref{defCond} for the definition of $\rho_n$).
Then $\rho$ is realized by a unique $\mu\in\mathcal{M}^*(\ddot\Gamma(\RR^d))$ if and only if the following hold.
\begin{eqnarray}
&\ &\hspace{-0.5cm} \tilde{L}_\rho(h^2)\geq 0,\,\, \forall h\in\mathscr{P},\label{psd}\\
&\ &\hspace{-0.5cm} \tilde{L}_{\rho}(\Phi_{\varphi, n}(\eta) \rangle h^2(\eta))\geq 0, \forall h\in\mathscr{P}, \forall \varphi\in\C_c^{+,\infty}(\RR^d),\forall\,n\in\NN,\label{psd-shift}
\end{eqnarray}
where $\Phi_{\varphi, n}(\eta):=\langle \varphi^{\otimes n},\eta^{\odot n}\rangle$.
\end{theorem}
Also, combining Lemma \ref{detCorr} with Corollary \ref{ThmSimpleConfig} we get:

\begin{theorem}\label{ThmSimpleConfig-Corr}\ \\
Let $\rho\in\mathcal{F}\left(\mathcal{R}\right)$ be a Stieltjes determining sequence s.t. the corresponding $\left(\frac{\rho_n}{n!}\right)_{n\in\NN}$ is almost-increasing (see \eqref{defCond} for the definition of $\rho_n$).
Then $\rho$ is realized by a unique $\mu\in\mathcal{M}^*(\Gamma(\RR^d))$ \underline{if and only if} the following inequalities hold
\begin{eqnarray*}
&\ & \tilde{L}_\rho(h^2)\geq 0,\,\, \forall h\in\mathscr{P},\\
&\ & \tilde{L}_\rho(\Phi_{\varphi, 1} h^2)\geq 0,\,\,\forall h\in\mathscr{P},\,\, \forall \varphi\in\C_c^{+,\infty}(\RR^d),\\
&\ & \tilde{L}_\rho(\Phi_{\varphi, 2} h^2)\geq 0,\,\,\forall h\in\mathscr{P},\,\, \forall \varphi\in\C_c^{+,\infty}(\RR^d),\\
&\ & \rho^{(2)}(diag(\Lambda\times\Lambda))=0,\,\,\,\forall\,\Lambda\in\mathcal{B}(\RR^d)\,\text{compact},
\end{eqnarray*}
where $diag(\Lambda\times\Lambda):=\{(\x,\x)\,|\, \x\in\Lambda\}$, $\Phi_{\varphi, 1}(\eta)=\langle \varphi,\eta\rangle$ and $\Phi_{\varphi, 2}(\eta)=\langle \varphi^{\otimes 2},\eta^{\odot 2}\rangle$.
\end{theorem}
Clearly the previous results also hold when $\rho=(\rho^{(j)})_{j\in\NN}\in\mathcal{F}\left(\mathcal{R}\right)$ fulfills \eqref{StieltjGen} and if the corresponding $\left(\frac{\rho_n}{n!}\right)_{n\in\NN}$ is almost-increasing (see \eqref{defCond} for the definition of $\rho_n$).

We can easily see that the conditions \eqref{psd} and \eqref{psd-shift} can be interpreted as that the sequence $\rho=(\rho^{(n)})_{n\in\NN_0}$ and all its shifted versions $_{\Phi_{\varphi, n}}\rho$ are positive semi-definite in the sense of Definition~\ref{PosSemiDef} by just using the following shift. 

\begin{definition}\label{SHIFT}\ \\
Given a sequence $\rho\in\mathcal{F}\left(\mathcal{R}\right)$ and a polynomial $P\in\mathscr{P}$ of the form $P(\eta) := \sum_{j=0}^{N}\langle p^{(j)},\eta^{\odot j}\rangle$, we define the sequence $_P\rho=\left((_P\rho)^{(n)}\right)_{n\in\NN_0}$ in $\mathcal{F}\left(\mathscr{D}'\right)$ as follows
$$\forall Q\in\mathscr{P},\quad \tilde{L}_{_{P}\rho}(Q):=\tilde{L}_\rho(PQ).
$$
\end{definition}

We intend now to understand more concretely the action of this shift. First of all, let us observe the following property whose proof is postponed to Appendix \ref{App1}.

\begin{lemma}\label{lemshifta}
For any $k, n\in\NN$ and any $\varphi\in\C_c^{+,\infty}(\RR^d)$, let $h^{(k)}\in\C_c^\infty(\RR^{{k}d})$ and $\Phi_{\varphi, n}(\eta)$ defined as above. Then we have:
$$\langle h^{(k)},(_{\Phi_{\varphi, n}}\rho)^{(k)}\rangle=\sum_{l=0}^{n\land k}\frac{n!k!}{(n-l)!l!(k-l)!}\left\langle \left[h^{(k)}\left(\varphi^{\otimes l}\hat{\otimes}1^{\otimes k-l}\right)\right]\hat{\otimes}\varphi^{\otimes n-l},\rho^{(n-l+k)}\right\rangle.$$
\end{lemma}

Let us take any polynomial $h(\eta):=\sum_{j=0}^\infty\frac{1}{j!}\langle h^{(j)}, \eta^{\otimes j}\rangle\in\mathscr{P}$ and let us denote by $H$ the sequence of coefficients of $h$, i.e. $H:=(h^{(j)})_{j=0}^\infty\in\F_f(\C_c^\infty)$. Then by using the definition and the properties of $K-$transform together with Lemma~\ref{lemshifta}, we get:
\begin{eqnarray*}
\tilde{L}_{\rho}(\Phi_{\varphi, n}(\eta)h^2(\eta))
&=&\tilde{L}_{_{\Phi_{\varphi, n}}\rho}(K(H\star H))\\
&=&\sum_{k=0}^\infty\frac{1}{k!}\sum_{l=0}^{n\land k}\frac{n!k!}{(n-l)!l!(k-l)!}\langle \left((H\star H)^{(k)}\varphi^{\otimes l}\hat{\otimes}1^{\otimes k-l}\right)\hat{\otimes}\varphi^{\otimes n-l},\rho^{(n-l+k)}\rangle,\\
&=&\sum_{i=0}^{\infty}\frac{1}{i!}\langle\sum_{l=0}^n\frac{n!}{(n-l)!l!} \left((H\star H)^{(i+l)}1^{\otimes i}\hat{\otimes}\varphi^{\otimes l}\right)\hat{\otimes}\varphi^{\otimes n-l},\rho^{(i+n)}\rangle,
\end{eqnarray*}
where in the last step we used the change of variables $i=k-l$.
Hence, the condition $$\tilde{L}_{\rho}(\Phi_{\varphi, n}(\eta)h^2(\eta))\geq 0,\quad \forall h\in \mathscr{P},$$ becomes:
\begin{equation}\label{newcond*}
\sum_{i=0}^{\infty}\frac{1}{i!}\langle\sum_{l=0}^n\frac{n!}{(n-l)!l!} \left((H\star H)^{(i+l)}\varphi^{\otimes l}\hat{\otimes}1^{\otimes i}\right)\hat{\otimes}\varphi^{\otimes n-l},\rho^{(i+n)}\rangle\geq 0, \forall H\in \F_f(\C_c^\infty).
\end{equation}
In the following we would like to rewrite this positive semi-definite condition in a new way which separates the action of $\varphi$ from the one of $H$ in order to make the condition more effective in some concrete cases of interest in applications (see e.g. Example \ref{PoissonProc}).

Let us first observe that the test function in the pairing can be rewritten as:
\begin{eqnarray}\label{rewriting-test-function-inside}
&\ &\sum_{l=0}^n\frac{n!}{(n-l)!l!} \left((H\star H)^{(i+l)}1^{\otimes i}\hat{\otimes}\varphi^{\otimes l}\right)\hat{\otimes}\varphi^{\otimes n-l}(y_1,\ldots, y_i, x_1,\ldots, x_n)=\\\nonumber
&=&(\tilde{H}_{x_1,\ldots, x_n}\star\tilde{H}_{x_1,\ldots, x_n})^{(i)}(y_1,\ldots, y_i)\varphi^{\otimes n}(x_1,\ldots, x_n)
\end{eqnarray}
where for any $s\in\NN$ we have set:
\begin{equation}\label{Htilde}
\tilde{H}_{x_1,\ldots, x_n}^{(s)}(y_1,\ldots, y_s):=\sum_{J\subset \{1,\ldots, n\}}H^{(s+|J|)}(y_1,\ldots, y_s, x_{J})
\end{equation} (see Appendix \ref{App1} for a proof of this rewriting).

This suggests that the condition 
$$\sum_{i=0}^{\infty}\frac{1}{i!}\langle (H\star H)^{(i)}\hat{\otimes}\varphi^{\otimes n}, \rho^{(i+n)} \rangle \geq 0, \forall {H}\in \F_f(\C_c^\infty)$$
may imply \eqref{newcond*}. This is the case whenever each $\rho^{(k)}$ is absolutely continuous with respect to $\sigma^{\otimes k}$ for some $\sigma\in\R(\RR^d)$. Indeed, we have the following weaker version of Theorem~\ref{ThmMultiplConfig-Corr} and \ref{ThmSimpleConfig-Corr}.
\begin{theorem}\label{Thm-Suff}
Let $\rho\in\mathcal{F}\left(\mathcal{R}\right)$ be a Stieltjes determining sequence s.t. each $\rho^{(k)}$ is absolutely continuous with respect to $\sigma^{\otimes k}$ for some $\sigma\in\R(\RR^d)$ and the corresponding $\left(\frac{\rho_n}{n!}\right)_{n\in\NN}$ is almost-increasing (see \eqref{defCond}). \\
\noindent(a) If the following holds for all $n\in\NN_0$, for $\sigma\mbox{-a.a. } x_1,\ldots, x_n \in \RR^d$ and for all $H\in \F_f(\C_c^\infty)$:
\begin{equation}\label{cond-Thm-Suff}
\sum_{i=0}^{\infty}\frac{1}{i!}\int_{\mathbb{R}^{id}} (H\star H)^{(i)}(y_1, \ldots , y_i)  \frac{d\rho^{(i+n)}}{d\sigma^{i+n}}(x_1, \ldots , x_n, y_1, \ldots , y_i) \sigma(dy_1) \ldots \sigma(dy_i)  \geq 0,
\end{equation}
then $\rho$ is realized by a unique $\mu\in\mathcal{M}^*(\ddot\Gamma(\RR^d))$.

\noindent (b) If $\sigma$ is non-atomic and \eqref{cond-Thm-Suff} holds for $n=0,1,2$, for $\sigma\mbox{-a.a. } x_1,\ldots, x_n \in \RR^d$ and for all $H\in \F_f(\C_c^\infty)$, then $\rho$ is realized by a unique $\mu\in\mathcal{M}^*(\Gamma(\RR^d))$ .
\end{theorem}
\proof
First let us observe that, since each $\rho^{(k)}$'s is absolutely continuous w.r.t. $\sigma^{\otimes k}$, \eqref{cond-Thm-Suff} for $n=0$ always holds by just applying the definition of the Radon-Nikodym derivative, i.e. 
$$\sum_{i=0}^{\infty}\frac{1}{i!}\langle (H\star H)^{(i)}, \rho^{(i)} \rangle \geq 0, \forall H\in \F_f(\C_c^\infty),$$
always holds and so \eqref{psd} holds.

Let $n\in\NN$. By applying the assumption \eqref{cond-Thm-Suff} for $H^{(i)}=(\tilde{H}_{x_1,\ldots, x_n})^{(i)}$, we get for $\sigma\mbox{-a.a. } x_1,\ldots, x_n \in \RR^d$ that the following holds
$$\sum_{i=0}^{\infty}\frac{1}{i!}\int_{\mathbb{R}^{id}} (\tilde{H}_{x_1,\ldots, x_n}\star \tilde{H}_{x_1,\ldots, x_n})^{(i)}(y_1, \ldots , y_i)  \frac{d\rho^{(i+n)}}{d\sigma^{i+n}}(x_1, \ldots , x_n, y_1, \ldots , y_i) \sigma(dy_1) \ldots \sigma(dy_i)  \geq 0,$$
where $\tilde{H}_{x_1,\ldots, x_n}$ is defined as in \eqref{Htilde}. Hence, for any $\varphi\in\C_c^{+,\infty}(\RR^d)$ multiplying the integrand by the non-negative function $\varphi^{\otimes n}$ and integrating w.r.t. $\sigma^{\otimes n}$ we obtain that:
 $$\sum_{i=0}^{\infty}\frac{1}{i!}\int_{\RR^{dn}}(\tilde{H}_{x_1,\ldots, x_n}\star\tilde{H}_{x_1,\ldots, x_n})^{(i)}(y_1,\ldots, y_i)\varphi^{\otimes n}(x_1,\ldots, x_n)\rho^{(i+n)}(dy_1,\ldots, dy_i, dx_1,\ldots, dx_n)\geq 0$$
By \eqref{rewriting-test-function-inside}, this exactly gives \eqref{newcond*} which is equivalent to \eqref{psd-shift} as
we have showed above. Then the conclusion (a) is a consequence of Theorem~\ref{ThmMultiplConfig-Corr}, while (b) of Theorem~\ref{ThmSimpleConfig-Corr} considering that $\rho^{(2)}(diag(\Lambda\times\Lambda))=0$ follows from the fact that $\rho^{(2)}$ is absolutely continuous w.r.t. $\sigma^{\otimes 2}$ and $\sigma$ is non-atomic.\\
\endproof

Let us present now an application of this theorem which illustrates a concrete use of this rewriting of the realizability conditions.

\begin{example}\label{PoissonProc}
Let $\sigma$ be a non-negative Radon measure on $\RR^d$. The Poisson measure $\mu_\sigma$ on $\mathscr{D}'(\RR^d)$ is defined as the unique measure on $\mathscr{D}'(\RR^d)$ which has $\sigma^{\otimes n}$ as correlation function $\rho^{(n)}$. Using Minlos theorem one can easily see that $\mu_\sigma$  is the unique measure on $\mathscr{D}'(\RR^d)$ with this property, since the infinite dimensional Laplace-Fourier transform associated to $\mu_\sigma$  is given by 
\[
\exp\left(\int_{\RR^d} (e^{\varphi(x)}-1) \sigma(dx)\right),
\]
which is continuous and positive semi-definite. The latter is a direct consequence of the positive semi-definiteness of the function $y \mapsto e^{z(e^{y}-1)}$ on $\mathbb{R}$ for $z \geq 0$, see e.g. Subsections~III.4.2 and 4.3 in \cite{GV68} and \cite{KMM78, Ka83}. It is less easy to see without an explicit construction that $\mu_\sigma$ is supported on the $\ddot\Gamma(\RR^d)$.
{However, by Theorem \ref{Thm-Suff} we can conclude that $\mu_\sigma$ is supported on $\ddot\Gamma(\RR^d)$ and even on  $\Gamma(\RR^d)$ whenever $\sigma$ is non-atomic.} Indeed,  as the $i-$th correlation function $\rho^{(i)}$ of the Poisson measure $\mu_\sigma$ is given by the product measure $\sigma^{\otimes i}$, the Radon-Nykodym derivative is just equal to $1$. Hence, condition \eqref{cond-Thm-Suff} is in this case independent of $n$, which means that the positive semi-definiteness of all the shifted sequences $(\rho^{(i+n)})_{i\in\NN}$ for all $n\in\NN$ follows just by the positive semi-definiteness of $\rho$. The latter is a direct consequence of the positive semi-definiteness of the Laplace-Fourier transform.
\end{example}

\section{Appendix}
For the convenience of the reader, we collect in this appendix some technical proofs {of auxiliary results used throughout this article}.

\subsection{{Proofs of some auxiliary results}}\label{App1}\ \\
Let us first give a proof of Lemma \ref{lemmino} which provides the construction of a total subset of $\C_c^\infty(\RR)$ closed under multiplication.\\

\emph{Proof of Lemma \ref{lemmino}}\ \\
For any $y,p\in\mathbb{Q} $ we have that
\begin{align*}
  &(\| f_{y,p,j}\|_{H_{k^{(n)}}})^2 \leq  \sum_{k=0}^{k_1^{(n)}} \int_{\mathbb{R}} \left| \frac{d^k}{dx^k}\varphi^j(x-y)e^{ipx}\right|^2  k_2^{(n)}(x) dx \\
  & \leq  \sum_{k=0}^{k_1^{(n)}} \int_{\mathbb{R}}\left( \sum_{n_0+\cdots+n_j=k} \frac{k!}{n_0!\cdots n_j!}\left| \frac{d^{n_0}}{dx^{n_0}}e^{ipx}\right| \prod_{i=1}^j \left|\frac{d^{n_i}}{dx^{n_i}}\varphi^j(x-y)\right|\right)^2  k_2^{(n)}(x) dx.
  \end{align*}
 W.l.o.g. we can always assume that the log-convex sequence $(d_l)_l$ is monotone increasing with $d_0=1$. Using the bound for derivative of $\varphi$ and the properties the sequence $(d_l)_l$ we get
 \begin{align*}
 (\| f_{y,p,j}\|_{H_{k^{(n)}}})^2& \leq  \sum_{k=0}^{k_1^{(n)}} \int_{\mathbb{R}}\left( \sum_{n_0+\cdots+n_j=k} \frac{k!}{n_0!\cdots n_j!} |p|^{n_0}\prod_{i=1}^jd_{n_i}\right)^2  k_2^{(n)}(x) dx\\
 & \leq  \sum_{k=0}^{k_1^{(n)}} \int_{\mathbb{R}}\left( \sum_{n_0+\cdots+n_j=k} \frac{k!}{n_0!\cdots n_j!} |p|^{n_0}d_{k-n_0}d_0\right)^2  k_2^{(n)}(x) dx \\
 & \leq d_{k_1{^{(n)}}}^2 \sum_{k=0}^{k_1^{(n)}} \int_{\mathbb{R}}\left( \sum_{n_0+\cdots+n_j=k} \frac{k!}{n_0!\cdots n_j!} |p|^{n_0}\right)^2  k_2^{(n)}(x) dx \\
 & \leq d_{k_1{^{(n)}}}^2 \sum_{k=0}^{k_1^{(n)}}  (j+|p|)^{2k}
 \int_{[-1,1]}     k_2^{(n)}(x+y) dx \\
 &\leq d_{k_1{^{(n)}}}^2 (j+|p|)^{k_1^{(n)}} k_1{^{(n)}}\sup_{x\in[-1,1]}|k_2^{(n)}(x+y) |
   \end{align*}
    
The totality of $E_0$ in $\mathscr{D}(\mathbb{R})$ is given by the second part of the proof of \cite[Lemma~2.5]{IKR14}.
\endproof

We prove now the form of the shift introduced in Definition~\ref{SHIFT} and stated in Lemma \ref{lemshifta}.

\noindent{\it Proof of Lemma \ref{lemshifta}}\\
We set $G:=(g^{(j)})_{j\in\NN}\in\F_f(\C^\infty_c)$ to be such that $g^{(j)}=0$ for all $j\neq n$ and $g^{(n)}=\varphi^{\otimes n}$ and $H:=(h^{(j)})_{j\in\NN}\in\F_f(\C^\infty_c)$ to be such that $h^{(j)}=0$ for all $j\neq k$.
\begin{eqnarray}
\langle h^{(k)},(_{\Phi_{\varphi, n}}\rho)^{(k)}\rangle\nonumber
&=&\tilde{L}_{_{\Phi_{\varphi, n}}\rho}(\langle h^{(k)},\eta^{\otimes k}\rangle)=\tilde{L}_{\rho}(\Phi_{\varphi, n}(\eta)\langle h^{(k)},\eta^{\otimes k}\rangle)\\ \nonumber
&=&\tilde{L}_{\rho}(n!KG(\eta)\cdot k!KH(\eta))\\
&=&\tilde{L}_{\rho}\left(\sum_{j=0}^\infty\frac{n!k!}{j!}\langle(G\star H)^{(j)}, \eta^{\otimes j}\rangle\right)
=\sum_{j=0}^\infty\frac{n!k!}{j!}\langle(G\star H)^{(j)}, \rho^{(j)}\rangle.\label{identitystar}
\end{eqnarray}
Now let us observe that for all $j\geq n\lor k$ we have:
\begin{eqnarray}
(G\star H)^{(j)}(x_J)\nonumber
&=&\sum_{\substack{J_1,J_2\subset J\\
                  J_1\cup J_2=J\\ |J_1|=n, |J_2|=k}}G^{(n)}(x_{J_1})H^{(k)}(x_{J_2})\\ \nonumber
&=&\sum_{\substack{I_1,I_2, I_3\subset J\\
                  I_1\cup I_2\cup I_3=J,
                 I_s \cap I_t=\emptyset,s\neq t\\
                |I_1|=j-k, |I_2|=n+k-j, |I_3|=j-n}}\varphi^{\otimes j-k}(x_{I_1})\varphi^{\otimes n+k-j}(x_{I_2})h^{(k)}(x_{I_2\cup I_3})\\ \nonumber
&=&\sum_{\substack{I_1,I_2, I_3\subset J\\
                  I_1\cup I_2\cup I_3=J,
                 I_s \cap I_t=\emptyset,s\neq t\\
                |I_1|=j-k, |I_2|=n+k-j, |I_3|=j-n}}\varphi^{\otimes j-k}\otimes\left((\varphi^{\otimes n+k-j}\otimes 1^{\otimes j-n})h^{(k)}\right)(x_{I_1\cup I_2\cup I_3})\\ 
&=&\frac{j!}{(j-k)!(n+k-j)!(j-n)!} \varphi^{\otimes j-k}\hat{\otimes}\left((\varphi^{\otimes n+k-j}\hat{\otimes} 1^{\otimes j-n})h^{(k)}\right)(x_J)\label{identity-intermediate}
\end{eqnarray}

Using \eqref{identity-intermediate} in \eqref{identitystar}, we get
\begin{eqnarray*}
\langle h^{(k)},(_{\Phi_{\varphi, n}}\rho)^{(k)}\rangle&=&\sum_{j=n\lor k}^{n+k}\frac{n!k!}{(j-k)!(n+k-j)!(j-n)!}\langle h^{(k)}\left(\varphi^{\otimes n+k-j}\hat{\otimes}1^{\otimes j-n}\right)\hat{\otimes}\varphi^{\otimes j-k},\rho^{(j)}\rangle \\
\end{eqnarray*}
The result follows using the change of variables $l=n+k-j$.
\endproof
In conclusion, we give the details for the rewriting of the realizability conditions introduced at the end of Section~3.\\

\emph{Proof of \eqref{rewriting-test-function-inside}.}
\begin{eqnarray*}
&\ &\sum_{l=0}^n\frac{n!}{(n-l)!l!} \left((H\star H)^{(i+l)}1^{\otimes i}\hat{\otimes}\varphi^{\otimes l}\right)\hat{\otimes}\varphi^{\otimes n-l}(y_1,\ldots, y_i, x_1,\ldots, x_n)=\\
&=&\sum_{l=0}^n\frac{n!}{(n-l)!l!} (H\star H)^{(i+l)}(y_1,\ldots, y_i, x_1,\ldots, x_l)\varphi^{\otimes l}(x_1,\ldots, x_l)\varphi^{\otimes n-l}(x_{l+1},\ldots, x_n)\\
&=&\!\!\!\!\!\sum_{J\subset\{1,\ldots, n\}}\sum_{\substack{Y_1,Y_2\subset \{1,\ldots, k\}\\
                  Y_1\cup Y_2=\{1,\ldots, k\}}}\sum_{\substack{J_1,J_2\subset J\\
                  J_1\cup J_2=J}} H^{(|Y_1|+|J_1|)}(y_{Y_1}, x_{J_1})H^{(|Y_2|+|J_2|)}(y_{Y_2}, x_{J_2})\varphi^{\otimes |J|}(x_{J})\varphi^{\otimes n-|J|}(x_{\{1,\ldots, n\}\setminus J})\\
&=&\!\!\sum_{\substack{Y_1,Y_2\subset \{1,\ldots, i\}\\
			   Y_1\cup Y_2=\{1,\ldots, i\}}}\!\!\sum_{J_1\subset \{1,\ldots, n\}}
			   H^{(|Y_1|+|J_1|)}(y_{Y_1}, x_{J_1})\!\!\sum_{J_2\subset \{1,\ldots, n\}}H^{(|Y_2|+|J_2|)}(y_{Y_2}, x_{J_2})\varphi^{\otimes n}(x_1,\ldots, x_n)\\
&=&(\tilde{H}_{x_1,\ldots, x_n}\star\tilde{H}_{x_1,\ldots, x_n})^{(i)}(y_1,\ldots, y_i)\varphi^{\otimes n}(x_1,\ldots, x_n)
\end{eqnarray*}

\subsection{Comparison with \cite{GIKM}}\label{Rem-SP}\  \\
A result similar to Theorem \ref{MainThm1} can be obtained by applying \cite[Theorem 1.1]{GIKM} for $V=\C_c^\infty(\RR^d)$ endowed with the projective topology $\tau_{proj}$ and $M=\mathcal{Q}(\mathscr{P}_{K})$, keeping in mind that $\mathscr{P}$ is algebraically isomorphic to the symmetric tensor algebra $S(V)$. In this way one gets that the conditions in \eqref{MainThmCond1-1} are necessary and sufficient to solve the KMP for $m$ under the assumption that $L_m$ is $\overline{\tau_{proj}}$-continuous, where $\overline{\tau_{proj}}$ is the natural extension of $\tau_{proj}$ to $\mathscr{P}$ considered in \cite{GIKM}. By \cite[Remark (9)-(10)]{GIKM}, it is clear that the $\overline{\tau_{proj}}$-continuity of $L_m$ implies that $m$ is a Stieltjes determining sequence. Note that this continuity assumption also forces the support of the representing measure to be contained in a compact subset of~$K$.

Also Theorem \ref{ThmSP} can be obtained by using \cite[Theorem 1.1]{GIKM} applied for $V=\C_c^\infty(\RR^d)$ endowed with the projective topology $\tau_{proj}$, but we need to take $M$ to be the quadratic module in $\mathscr{P}$ generated by the family $\{\Phi_\psi : \psi\in\C_c^{+,\infty}(\RR^d)\}\cup\{ \Upsilon_\varphi: \varphi\in\C_c^{+,\infty}(\RR^d), \left\|\varphi\right\|_\infty\leq 1\}$, where $\Phi_\psi$ and $\Upsilon_\varphi$ are the ones defined in Theorem \ref{ThmSP}. By the first part of the sufficiency proof and the Cauchy-Schwartz inequality, we get that:
\begin{equation}\label{rel}
L_m(\langle \varphi, \eta\rangle^{n})\leq L_m(1),\forall n\in\NN, \forall \varphi\in\C_c^{+,\infty}(\RR^d) \text{ with } \left\|\varphi\right\|_\infty\leq 1.
\end{equation}
We will now show that \eqref{rel} implies that $L_m$ is $\overline{\tau_{proj}}-$continuous.

Using the linearity of $L_m$, we can easily derive from \eqref{rel} that 
\begin{equation}\label{rel1}L_m(\langle \varphi, \eta \rangle^n) \leq L_m(1) \|\varphi\|_\infty^n,\,\forall n\in\NN, \forall \varphi\in\C_c^{+,\infty}(\RR^d).\end{equation}
Then for all $\varphi  \in\C_c^{\infty}(\RR^d)$ we have $\varphi=\varphi_+-\varphi_-$ with $\varphi_+(x):=\max\{\varphi(x), 0\}$, $\varphi_-(x):=-\min\{\varphi(x), 0\}$, and so:
\[
\left| L_m(\langle \varphi, \eta \rangle^n) \right|  \leq \sum_{k=0}^n
 \binom{n}{k} \left| L_m(\langle \varphi_-, \eta \rangle^k \langle \varphi_+, \eta \rangle^{n-k})\right| \]
By using Cauchy-Schwarz inequality in the previous inequality and then \eqref{rel1}, we get:
\be\label{rel2}
\left| L_m(\langle \varphi, \eta \rangle^n) \right|  \leq\sum_{k=0}^n
 \binom{n}{k} L_m(1) \|\varphi_-\|_\infty^k \|\varphi_+\|_\infty^{n-k}\leq 2^n L_m(1) 
  \|\varphi\|_\infty^n.
  \ee
This gives in turn that for each $n$ the following multilinear form on $\left(\C_c^\infty(\RR^d)\right)^n$ endowed with the product topology induced by $( \C_c^\infty(\RR^d), \|\cdot\|_\infty)$ is continuous: 
\[
(\varphi_1, \ldots , \varphi_n) \mapsto L_m(\prod_{i=1}^n\langle \varphi_i, \eta \rangle).\]
Indeed, using the polarization identity and then \eqref{rel2}, we have that:
\begin{eqnarray*}
\left| L_m(\prod_{i=1}^n\langle \varphi_i, \eta \rangle)
 \right| &=& \left|\frac{1}{2^n n!}\sum_{\varepsilon_1, \ldots , \varepsilon_n \in \{-1,1\}} \prod_{i=1}^n \varepsilon_i L_m\left( \langle \sum_{i=1}^n \varepsilon_i \varphi_i , \eta \rangle^n \right) \right|\\
&\leq& \frac{1}{n!} \sup_{\varepsilon_1, \ldots , \varepsilon_n \in \{-1,1\}} \left| L_m\left( \langle \sum_{i=1}^n \varepsilon_i \varphi_i , \eta \rangle^n \right) \right|\\
& \stackrel{\eqref{rel2}}{\leq}& \frac{L_m(1)2^n}{n!} \sup_{\varepsilon_1, \ldots , \varepsilon_n \in \{-1,1\}}  \left\| \sum_{i=1}^n \varepsilon_i \varphi_i\right\|_\infty^n,
\end{eqnarray*}
and therefore
\[
\sup_{\stackrel{\varphi_i\in\C_c^\infty(\RR^d)}{ \| \varphi_i\| \leq 1}}\left| L_m(\prod_{i=1}^n\langle \varphi_i, \eta \rangle)
 \right| \leq L_m(1) \frac{2^n n^n}{n!}\leq L_m(1) \frac{2^n e^n}{\sqrt{2\pi n} } \leq L_m(1) (2e)^n,
\]
where we used first that by triangle inequality we have $\left\| \sum_{i=1}^n \varepsilon_i \varphi_i\right\|_\infty \leq n$ and then the Stirling formula.

Hence, the multilinear mapping
\[
\begin{array}{ccl}
( \C_c^\infty(\RR^d), \|\cdot\|_\infty)\times\cdots\times( \C_c^\infty(\RR^d), \|\cdot\|_\infty)&\to&\RR\\
(\varphi_1, \ldots , \varphi_n) &\mapsto& L_m(\prod_{i=1}^n\langle \varphi_i, \eta\rangle).
\end{array}
\]
is continuous and so, by the universal property of the projective tensor product (see \cite[Proposition~43.4]{Tre67}), there exists a unique continuous functional $\tilde{L}_m$ on $\C_c^\infty(\RR^d)^{\otimes n}$ endowed with the projective tensor norm $\|\cdot\|_{\infty}^{\otimes n}$ such that
\[
L_m(\prod_{i=1}^n\langle \varphi_i, \eta \rangle) = \tilde{L}_m(\langle\varphi_1 \otimes \ldots \otimes \varphi_n, \eta^{\otimes n}  \rangle).
\]
In particular, this holds for $S(\C_c^\infty(\RR^d))_n\subset \C_c^\infty(\RR^d)^{\otimes n}$ and so by \cite[Proposition 43.12-b)]{Tre67} we have that
for all $f^{(n)}\in S(\C_c^\infty(\RR^d))_n$ 
\[
\left|  \tilde{L}_m(\langle f^{(n)}, \eta^{\otimes n}\rangle) \right| \leq L_m(1) (2e)^n \|f^{(n)}\|_{\infty}^{\otimes n}.
\]
Defining $\rho(\varphi) := 2 e \| \varphi \|_\infty,\, \forall\varphi\in\C_c^\infty(\RR^d)$, we have that
for all $f^{(n)}\in S(\C_c^\infty(\RR^d))_n$ 
\[
\left|  \tilde{L}(\langle f^{(n)}, \eta^{\otimes n}\rangle) \right| \leq L_m(1) {\rho^{\otimes n}}(f^{(n)})=L_m(1) \overline{\rho}_n(f^{(n)}),
\]
where $\rho^{\otimes n}$ is the projective tensor norm on $\C_c^\infty(\RR^d)^{\otimes n}$ induced by $\rho$ and $\overline{\rho}_n$ the quotient norm on $S(\C_c^\infty(\RR^d))_n$ induced by $\rho^{\otimes n}$.

Since every polynomial $p\in\mathscr{P}$ is of the form $p(\eta)=\sum_{n=1}^N\langle f^{(n)}, \eta^{\otimes n}\rangle$, we have
\[
\left|  \tilde{L}( p(\eta)) \right| =\left|  \tilde{L}( \sum_{n=1}^N\langle f^{(n)}, \eta^{\otimes n}\rangle) \right| \leq L_m(1) \sum_{n=1}^N\overline{\rho}_n(f^{(n)})=L_m(1)\overline{\rho}(p),
\]
where $\overline{\rho}:=\sum_{n=1}^N\overline{\rho}_n$ is exactly the extension used in \cite{GIKM} (keeping in mind that $S(V)\approx \mathscr{P}$).

Hence, we have showed that $L_m$ is $\overline{\rho}$-continuous and so $\overline{\tau_{proj}}-$continuous (see \cite[Chapter I, Section~3.10]{B86}  and \cite[Section 5]{GIKM} for more details on $\tau_{proj}$ and $\overline{\tau_{proj}}$, respectively). Then, by \cite[Theorem~1.1]{GIKM} \eqref{(a)}, \eqref{(b)} and \eqref{(c)} are equivalent to the existence of a unique representing measure for $L_m$ whose support is contained in $\mathcal{R}_{sub}(\RR^d)$, i.e. Theorem~\ref{ThmSP} holds.

\section*{Acknowledgments}
The work of M. Infusino and T. Kuna was partially supported by the EPSRC Research Grant EP/H022767/1, the Zukunftskolleg Mentorship Programme 2016 and the Young Scholar Fund of University of Konstanz 83942718. We are also indebted to the Baden-W\"urttemberg Stiftung for the financial support to this work by the Eliteprogramm f\"ur Postdocs 1.16101.17. The authors would like to thank {Abdelmalek Abdesselam}, Eugene Lytvynov and Maria Joao Oliveira for the helpful discussions.

\def\cprime{$'$}

\end{document}